\documentclass[11pt]{article}
\usepackage[a4paper,left=25mm,right=25mm,top=2cm,bottom=2cm]{geometry}
\usepackage[T1]{fontenc}
\usepackage{lmodern}
\usepackage{xcolor,microtype,mathtools,amssymb,amsthm,mathrsfs}
\usepackage{xurl}
\usepackage[colorlinks=true,linkcolor=blue!55!black,
 citecolor=blue!55!black,urlcolor=blue!55!black,
 bookmarksnumbered=true,bookmarksopen=true,
 bookmarksopenlevel=2,
 pdfpagemode=UseOutlines,pdfdisplaydoctitle=true]{hyperref}

\usepackage{bookmark}
\numberwithin{equation}{section}
\newtheorem{theorem}{Theorem}[section]
\newtheorem{proposition}[theorem]{Proposition}

\newtheorem{corollary}[theorem]{Corollary}
\theoremstyle{definition}

\theoremstyle{remark}

\hypersetup{pdftitle={Dimension-free weak-type (1,1) bounds for coordinatewise maximal Riesz truncations}}
\title{\Large\bf Dimension-free weak-type $(1,1)$ bounds for coordinatewise maximal Riesz truncations}
\author{Sheng-Chen Mao}
\date{}
\begin{document}
\maketitle
\vspace{-1.8em}
\noindent\textbf{Abstract}. We establish the dimension-free weak-type $(1,1)$ estimate for the maximal   truncations of the first-order Riesz transforms on $\mathbb R^d$, with the separate truncation radius in each coordinate. To prove this, we combine a fractional obstacle decomposition with Mellin frequency analysis to control logarithmically smoothed truncations. A coupling argument then reduces the smoothing error to annular integrals, which can be bounded by using angular cancellation and thin-slab geometry.  Moreover, we also show that both maximal operators, with common or independent radii, admit $L^p$ norms comparable to $\max\{p,p/(p-1)\}$ with absolute constants for every $p\in(1,\infty)$.



\section{Introduction}
Let $d\ge1$ be an integer, and
\begin{equation}\label{rdef}
 \begin{aligned}
 c_d&=\frac{\Gamma((d+1)/2)}{\pi^{(d+1)/2}},
 & K(z)&=c_d\frac{z}{|z|^{d+1}}\quad(z\ne0),\\
 R_{i,t}f(x)&=\int_{|x-y|>t}K_i(x-y)f(y)\,dy,
 & &1\le i\le d,\ t>0 .
 \end{aligned}
\end{equation}
For $f\in L^p(\mathbb R^d)$, $1\le p<\infty$, the integral in \eqref{rdef} converges absolutely at each $t>0$. Write $\mathbf R_t f=(R_{i,t}f)_{i=1}^d$ and $\mathbf Rf=(R_if)_{i=1}^d$ for the vectors of truncations and principal values. With our Fourier convention as in Subsection  \ref{not}, $R_i$ has the multiplier $-i\xi_i/|\xi|$. The principal values exist almost everywhere for $f\in L^1$. See \cite[Chapter~II]{SteinSI70} for more details. We define the common-radius and independent-radius maximal operators by
\begin{equation}\label{defB}
 \mathcal A_df(x)=\sup_{t>0}|\mathbf R_t f(x)|_{\ell^2},
 \qquad
 \mathcal B_df(x)=\left(\sum_{i=1}^d\sup_{t>0}|R_{i,t}f(x)|^2\right)^{1/2}.
\end{equation}
By the principal value limit and the definitions in \eqref{defB}, we see that $|\mathbf R f|_{\ell^2}\le\mathcal A_df\le\mathcal B_df$ almost everywhere. Dominated convergence in \eqref{rdef} yields continuity in $t>0$ for $f\in L^p$. Hence one may take both suprema in \eqref{defB} over positive rational radii.

\subsection{Main result}
Our main result of this article is the following
\begin{theorem}\label{main}
There is an absolute constant $C$ such that, for every $d\ge1$, every  $f\in L^1(\mathbb R^d)$, and every $\lambda>0$,
\begin{equation}\label{w11}
 \lambda\big|\{x\in\mathbb R^d:\mathcal B_df(x)>\lambda\}\big|
 \le C\|f\|_1.
\end{equation}
\end{theorem}

By pointwise domination, the bound \eqref{w11} applies to $\mathcal A_d$ and each scalar maximal truncation. In Corollary~\ref{mweak} we extend the estimate to finite complex measures. The corresponding estimates with independent heights hold for conjugate Poisson integrals and for tangential derivatives of Poisson extensions by Corollaries~\ref{cweak} and~\ref{pgrad}, respectively. For $1<p<\infty$, we show in Proposition~\ref{pgrow} that the norms of $\mathcal A_d$ and $\mathcal B_d$ on $L^p$ are comparable to $p^*=\max\{p,p/(p-1)\}$ with absolute constants.

\subsection{Previous results and the endpoint question}

In 1983, Stein \cite{Stein83}  established dimension-free strong bounds for $\mathbf R$, and in \cite[p.~203, Problem~(b)]{Stein86} he asked whether the weak-type constants for the Riesz transforms are also independent of dimension. For a single principal value transform, Janakiraman \cite{Jan04} obtained a logarithmic bound. Spector and Stockdale \cite[Theorems~1.1 and~1.3]{SS21} subsequently obtained the same order by a different argument and together with Ouyang, they \cite[Theorem~1.1]{OSS26} lately deduced the weak-type constant for the principal value vector is at most $2$ for  real-valued functions. In the scalar case, Ba\~nuelos, Hormozi, and Volberg \cite[Theorem~1.1]{BHV26a} identified the exact weak-type $(1,1)$ norm with the Hilbert transform constant of Davis \cite{Dav74}. The scalar weak-type $(p,p)$ constants also agree with those of the Hilbert transform on real-valued functions when $1<p\le2$, by \cite[Theorem~1.1]{BHV26b}.

For   maximal truncations, classical singular integral theory ensures weak type $(1,1)$ in each dimension (cf. e.g. \cite[Corollary~3.8.2]{Gra24}). In the strong range $p\in(1,\infty)$, the vector variation bound of Gillespie and Torrea \cite[Theorem~B]{GT04} implies a dimension-free estimate for the common-radius operator $\mathcal A_d$. Kucharski and Wr\'obel \cite[Theorem~1.3 and Corollary~1.2]{KW23} proved this strong estimate directly and established a dimension-free $L^2$ bound for $\mathcal B_d$. Applying the general comparison in \cite[Theorem~1.1]{KWZ26} to the coordinate harmonics, we have
\begin{equation}\label{pcomp}
 \|\mathcal B_df\|_p\le C_p\|\mathbf R f\|_{L^p(\ell^2)},
 \qquad 1<p<\infty ,
\end{equation}
where $C_p$ is independent of $d$. One may take $C_p=O(p^3)$ as $p\to\infty$ and $C_p=O((p-1)^{-3})$ as $p\downarrow1$. Together with the vector Riesz bound \cite{BW95}, these estimates imply $\|\mathcal B_d\|_{L^p\to L^p}=O((p^*)^4)$ in the two limits. Kucharski, Wr\'obel, and Zienkiewicz \cite[p.~630, Remark~1 and (1-6)]{KWZ26} asked whether one could derive a bound in terms of $\|f\|_p$ directly, without the comparison \eqref{pcomp} with $\mathbf Rf$. Our proof uses the weak bound \eqref{w11} and a separate estimate for $p\ge3$ to obtain the sharp order $p^*$ for $\mathcal B_d$ in \eqref{pnorm}, which particularly provides an affirmative answer to their question.

For $d\ge2$, the positive radial factorization of Liu, Melentijevi\'c, 
and Zhu \cite[Theorems~2--4]{LMZ24}, together with the vector Riesz 
bound \cite{BW95}, yields $\|\mathcal A_d\|_{L^p\to L^p}=O(p)$ as 
$p\to\infty$. We establish the same growth rate for the larger 
operator $\mathcal B_d$ in \eqref{pnorm}. For higher-order transforms, 
comparison estimates for odd and even kernels were developed in 
\cite{MOPV10,MOV11}. Kucharski, Wr\'obel, and Zienkiewicz 
\cite[Theorem~1.1]{KWZ26} established dimension-free comparisons 
for higher-order Riesz transforms, with a more direct proof 
subsequently obtained by Kucharski and Kwa\'snicki 
\cite[Theorem~1.1]{KK26} using generalized spherical averages.

Recall that a bound for the principal value vector concerns the limit of $\mathbf R_t f$ as $t\downarrow0$. The common-radius endpoint discussed by Guti\'errez and Manfredi \cite[Section~5.7]{GM26} requires uniform control of $\mathbf R_t f$ over all $t>0$. In \eqref{w11} we obtain this control in weak $L^1$ for the larger operator $\mathcal B_d$, which allows a different measurable radius $t_i(x)$ in each coordinate.

For the endpoint, the choice of norm on the right also matters. Mateu, Orobitg, P\'erez, and Verdera \cite[Theorem~1]{MOPV10} characterized the $L^2$ comparison for smooth homogeneous odd kernels by spherical harmonic factorization. In the scalar Riesz case, Mateu and Verdera \cite[Theorem~2]{MV06} showed that no finite constant can bound $\|\sup_{t>0}|R_{i,t}f|\|_{1,\infty}$ by $\|R_if\|_1$ for all admissible $f$. Cao, Liu, Si, and Yan \cite[Theorem~1.9]{CLSY26} constructed examples in the real Hardy space and proved the comparison over dyadic radii with a dimension-dependent constant. Our bound \eqref{w11} utilizes the original data norm $\|f\|_1$ and applies to all positive radii.

\subsection{Notation and outline of the proofs} \label{not} 

Throughout the paper, the letters $C$ and $c$ denote positive absolute constants that may change from one occurrence to the next. We write $a\lesssim b$ if $a\le Cb$, and $a\asymp b$ if both $a\lesssim b$ and $b\lesssim a$. In qualitative existence and approximation arguments, constants may depend on the data and the dimension. We use the notation
\begin{equation*}
 L_d:=\log(d+2),\quad
 |D|:=(-\Delta)^{1/2},\quad H_t:=e^{t\Delta},\quad P_s:=e^{-s|D|}.
\end{equation*}
$H_t$ and $P_s$ denote both the semigroups and their convolution kernels. We write $H_*f=\sup_{t>0}|H_tf|$ and $P_*f=\sup_{s>0}|P_sf|$ for the associated scalar maximal functions. The symbol $H^s$ denotes the inhomogeneous $L^2$ Sobolev space.
For a measurable function $f$, we use the weak $L^1$ quasi-norm
\begin{equation*}
 \|f\|_{1,\infty}:=\|f\|_{L^{1,\infty}}:=\sup_{\lambda>0}\lambda|\{|f|>\lambda\}|.
\end{equation*}
We take the Fourier transform with phase $e^{-ix\cdot\xi}$ and inverse factor $(2\pi)^{-d}$. In the logarithmic radius $u$, the corresponding phase and inverse factor are $e^{-iu\omega}$ and $(2\pi)^{-1}$. We write $*_u$ for convolution in this variable and $b(D_u)$ for the Fourier multiplier with symbol $b$.

We equip $\mathbb C^d$ with the Euclidean norm and denote the corresponding space of vector-valued functions by $L^p(\ell^2)$. For $V=(V_i(u,x))_{i=1}^d$, write
\begin{equation}\label{nmax}
 \mathcal N V(x):=\left(\sum_i\sup_{u}|V_i(u,x)|^2\right)^{1/2}.
\end{equation}
Here $u$ ranges over the parameter domain of $V$. We denote the indicator of a set $E$ by $\mathbf1_E$ and use $\alpha:=(d+1)/2$ throughout.

This article is organized as follows. In Sections~\ref{s2} and~\ref{s3}, we prove scalar maximal bounds for positive extensions and an $L^2$ estimate for spherical first moments applied to densities of bounded variation in the radius. The former bounds enter the frequency analysis. The latter estimate is used when we change the center of the kernel. In Section~\ref{s4}, we control the contact trace and the energy of the Poisson extension. We then decompose in Mellin frequency in Section~\ref{s5}. At low frequencies we use the conjugate Poisson estimate \eqref{csoft} from Proposition~\ref{csf}. At intermediate frequencies we use the estimate for $W_s^{(m)}$ in Proposition~\ref{ascal}, with $m$ chosen for each band. At high frequencies the angle of the complex heat time compensates for the decay of the Gamma multiplier in \eqref{bray}.

To recover the hard cutoff we construct, in Section~\ref{s6}, the coupling $\Pi$ in \eqref{cdpi} between $g(y)\,dy$ and $h(z_*)\,dz_*$. The annular estimate \eqref{csurv} controls the tail of its conditional measures. In Section~\ref{s7}, we use the spherical estimate \eqref{bmest} to change the kernel center on pairs with small relative separation. On the complementary region we compare a hard truncation with two shifted smooth truncations. The resulting error $E_w$ is defined in \eqref{perr}. The bound \eqref{ppair} for $E_w^2$ follows by estimating the products of the opposite-sign contributions in each coordinate. When the source points are close, their angular separation contributes a quadratic factor. When they are separated, the two shell conditions restrict the evaluation point to a thin slab about their perpendicular bisector. Applying \eqref{csurv} to the annular densities combines the two contributions. With $w=cd^{-7/6}$, both the smoothed bound \eqref{bbest} and the error estimate \eqref{pthin} are uniform in $d$.

In Section~\ref{s8}, we determine the dependence of the strong bounds on $p$. For large $p$ we use a heat square function and a discrete strip representation of the conjugate Poisson integrals. Harmonic extension and boundary projection go back to Gundy and Varopoulos \cite{GV79} for the principal value transforms. In the strip representation each component can use a different stopping level, which corresponds to the independent choice of height.

\section{Scalar maximal estimates for positive extensions}\label{s2}

For $m>0$, let $W_s^{(m)}$ be the probability kernel
\begin{equation}\label{aker}
 W_s^{(m)}(x):=\frac{\Gamma((d+m)/2)}{\pi^{d/2}\Gamma(m/2)}
                   \frac{s^m}{(|x|^2+s^2)^{(d+m)/2}},\qquad s>0 .
\end{equation}
At $m=1$ we recover the Poisson kernel $P_s$, whose maximal operator satisfies a classical dimension-free weak $(1,1)$ estimate (see e.g. \cite[Chapter~II, Section~2, Lemma~1]{SteinLP70}). We use the dimension-free Poisson maximal inequality (see \cite[Chapter~III, Section~3]{SteinLP70})
\begin{equation}\label{apois}
 \left\|\sup_{s>0}|P_sq|\right\|_2\leq C\|q\|_2,
\end{equation}
and the weak $(1,1)$ estimate of Spector and Stockdale for the heat maximal operator (see \cite[Theorem~1.2]{SS26}),
\begin{equation}\label{aheat}
 \|H_*q\|_{L^{1,\infty}}\le C L_d\|q\|_1 .
\end{equation}
For general $m$, we express $W_s^{(m)}$ as an average of heat kernels against a gamma density. This representation permits a second bound in which $m$ replaces $d$ inside the logarithm.

\begin{proposition}\label{ascal}
For all $d\ge1$, $m>0$, and finite complex Borel measures $q$ on $\mathbb R^d$,
\begin{equation}\label{amax}
 \left\|\sup_{s>0}|W_s^{(m)}*q|\right\|_{L^{1,\infty}}
    \le C\min\{1+\log(1+m),L_d\}|q|(\mathbb R^d).
\end{equation}
\end{proposition}

\begin{proof}
Put $b=m/2$ and $v=\log(s^2/4)$. For a finite measure $\eta$, write $F_\eta(v)=H_{e^v}\eta$. Integration of the Gaussian against the gamma density yields
\begin{equation}\label{subor}
 W_s^{(m)}*\eta
 =\frac{1}{\Gamma(b)}\int_0^\infty H_{e^v/r}\eta\,r^{b-1}e^{-r}\,dr
 =(L_b*F_\eta)(v),\qquad L_b(w)=\frac{e^{bw-e^w}}{\Gamma(b)} .
\end{equation}
Since $L_b$ is a probability density in $v$, \eqref{subor} bounds $\sup_s|W_s^{(m)}*q|$ by $H_*|q|$ for $q\in L^1$. We thus obtain the factor $L_d$ from \eqref{aheat} and extend this estimate to finite measures below. For $d=1,2$ this already proves \eqref{amax}, so it remains to establish the bound with $1+\log(1+m)$ when $d\ge3$.

Write $D_m=1+\log(1+m)$ and begin with a finite positive atomic measure $\nu$ of mass $M$. For $\kappa>0$, Spector and Stockdale \cite[Proposition~2.1]{SS26} construct an open set $O$ and a nonnegative potential $U$ vanishing on $O^c$ such that $\kappa|O|=M$ and $\sigma:=\nu-\mu=-\Delta U$, where $\mu=\kappa\mathbf1_O\,dx$. For $G(\tau,x)=i\tau\widehat F_\sigma(\tau,x)$, their frequency estimate \cite[Theorem~1.4, (1.18)]{SS26} reads
\begin{equation}\label{ssloc}
 \int_{O^c}\int_{|\tau|\le B}|G(\tau,x)|^2\,d\tau\,dx
 \le C B\kappa M,\qquad B\ge1 .
\end{equation}
Let $e(v)=e^{-v}\mathbf1_{\{v>0\}}$. On $O^c$,
\begin{equation*}
 (e*F_\sigma)(v)=e^{-v}(U-H_{e^v}U)\le0 .
\end{equation*}
Indeed, $(e*F_\sigma)(v)$ is the average of $H_t\sigma$ over $0<t<e^v$. Moreover, we have $F_\sigma(\cdot,x)\in L^1(\mathbb R)$ for almost every $x\in O^c$  from \cite[Section~4, proof of (1.18)]{SS26}. For such $x$, we use $F_\mu\le\kappa$ and $L_b\ge0$ in \eqref{subor} to obtain
\begin{equation*}
 W_s^{(m)}*\nu\le\kappa+L_b*(F_\sigma-e*F_\sigma)(v)
 \qquad(x\in O^c).
\end{equation*}
Since $1-\widehat e(\tau)=i\tau/(1+i\tau)$, the Fourier transform in $v$ of the last convolution is $\widehat L_b(\tau)G(\tau,x)/(1+i\tau)$.

Substituting $r=e^w$ in the Fourier transform of $L_b$ from \eqref{subor}, we obtain $\widehat L_b(\tau)=\Gamma(b-i\tau)/\Gamma(b)$. For $b>0$, Euler's product formula (cf. e.g. \cite[formula~5.8.3]{NIST}) reads
\begin{equation}\label{agam}
 \frac{|\Gamma(b+i\tau)|}{\Gamma(b)}
 =\prod_{j\ge0}\left(1+\frac{\tau^2}{(b+j)^2}\right)^{-1/2}.
\end{equation}
For the factors with $b+j\ge\max(b,1,|\tau|)$, we use $\log(1+y)\ge y/2$ on $[0,1]$. All other factors in \eqref{agam} are at most one. Since the sum of $(b+j)^{-2}$ over the indicated indices is at least $c/(1+b+|\tau|)$, we obtain
\begin{equation*}
 |\widehat L_b(\tau)|=\frac{|\Gamma(b+i\tau)|}{\Gamma(b)}
 \le\exp[-c\min\{\tau^2/(1+m),|\tau|\}] .
\end{equation*}
On a dyadic band $I_B=\{B<|\tau|\le2B\}$, $B=2^j\ge1$, Cauchy--Schwarz in $\tau$ and \eqref{ssloc} imply
\begin{equation*}
 \left\|\int_{I_B}\frac{|\widehat L_b G|}{|1+i\tau|}\,d\tau\right\|_{L^2(O^c)}
 \le C B^{-1/2}\sup_{I_B}|\widehat L_b|\,(B\kappa M)^{1/2}.
\end{equation*}
The same argument on $|\tau|\le1$ yields $C\sqrt{\kappa M}$. There are $O(D_m)$ dyadic bands up to $B=\sqrt{1+m}$, and the exponential bound for $\widehat L_b$ makes $\sum_{B>\sqrt{1+m}}\sup_{I_B}|\widehat L_b|$ uniformly bounded. Summing over the bands and inverting the Fourier transform, we obtain
\begin{equation*}
 \left\|\sup_{v\in\mathbb R}|L_b*(F_\sigma-e*F_\sigma)(v)|\right\|_{L^2(O^c)}
 \le C D_m\sqrt{\kappa M}.
\end{equation*}
The band estimates ensure absolute integrability of the Fourier transform for almost every $x\in O^c$. At such points the inverse transform agrees everywhere with the continuous convolution in $v$. Taking $\kappa=\lambda/(2D_m)$, we bound the measure of the $\lambda$-superlevel set by $CD_mM/\lambda$, using $|O|=M/\kappa$ on $O$ and Chebyshev's inequality on $O^c$.

To pass from atoms to a compactly supported finite positive measure, partition its support by a cubical mesh and place the mass of each cube at its center. For each $s>0$, uniform continuity of $W_s^{(m)}$ ensures uniform convergence of the convolutions as the mesh size tends to zero. Fatou's lemma transfers the estimate with $D_m$ to maxima over finitely many radii. For the estimate with $L_d$, use the same argument with mollifications of the measure and the bound for $L^1$ functions. Exhausting $\mathbb R^d$ by compact sets and the positive rationals by finite sets now proves both estimates for every finite positive measure, since each convolution is continuous in $s>0$. The complex case follows by domination by the convolution with $|q|$.
\end{proof}

\section{Spherical averages and spatial frequency decomposition}\label{s3}

After the change of kernel center in Section~\ref{s7}, the density in a spherical integral will vary with the radius. We estimate the tails of this integral assuming a pointwise bound on the total variation of the density. Let $\sigma_{d-1}$ be normalized surface measure on $S^{d-1}$, and put $\alpha_d=c_d|S^{d-1}|\asymp\sqrt d$. We denote the spherical average of $a$ by $S_ra=\int_{S^{d-1}}a(\cdot-r\theta)\,d\sigma_{d-1}(\theta)$. For $d\ge2$, we also use the vector of first spherical moments
\begin{equation}\label{shell}
 D_r a(x):=\alpha_d\int_{S^{d-1}}\theta\,a(x-r\theta)\,d\sigma_{d-1}(\theta)
 \qquad(r>0).
\end{equation}
For $a\in L^2$ and $r>0$, we interpret $D_ra$ as an element of $L^2(dx;\ell^2)$. If $a_r$ is jointly measurable in $(r,x)$, integration in \eqref{shell} defines a jointly measurable representative wherever the integral converges absolutely. By translation invariance and Fubini's theorem, changing $a_r$ on a spatial null set for each $r$ affects this representative only on a null set in $(r,x)$.

The multiplier of $D_r$ changes scale when $r|\xi|\asymp\sqrt d$. To see this, let $\rho_j$ denote the density of the first coordinate under normalized surface measure on $S^{j-1}$, $j\ge2$, and write $\Phi_j=\widehat\rho_j$. For $\xi\ne0$ we have
\begin{equation}\label{smult}
 \widehat D_r(\xi)=-i\frac{\alpha_d}{d}\frac{\xi}{|\xi|}(r|\xi|)\Phi_{d+2}(r|\xi|).
\end{equation}
Indeed, by rotation invariance it suffices to take $\xi=|\xi|e_1$. All transverse components vanish by symmetry, and integration by parts yields \eqref{smult} because $\rho'_{d+2}(t)=-dt\rho_d(t)$, where $\rho_d(t)$ is proportional to $(1-t^2)_+^{(d-3)/2}$. To estimate the Fourier transform in \eqref{smult}, write $\rho=\rho_{d+2}$. For $d\ge4$, we have $\rho(t)=C_\gamma(1-t^2)_+^\gamma$ with $\gamma=(d-1)/2>1$. Evaluating the beta integrals, we find
\begin{align*}
 \|\rho''\|_1
 &\le 2\gamma\int_{-1}^1\frac{\rho(t)}{1-t^2}\,dt+
  4\gamma(\gamma-1)\int_{-1}^1\frac{t^2\rho(t)}{(1-t^2)^2}\,dt\\
 &=4(\gamma+1/2)=2d.
\end{align*}
Since $\rho$ and $\rho'$ vanish at both endpoints, $|\Phi_{d+2}(s)|\le\min\{1,2d/s^2\}$. The same estimate holds for $d=3$, since the distributional second derivative of $\rho(t)=3(1-t^2)_+/4$ has total variation $6$. For $d=2$, one integration by parts reduces the estimate to the Fourier transform of $\rho'$. When $s\ge2$, the $L^1$ norm of $\rho'$ on the two endpoint intervals of length $s^{-1}$ is $O(s^{-1/2})$. On the remaining interval, the boundary values of $\rho'$ and the integral of $|\rho''|$ are $O(s^{1/2})$. A second integration by parts on that interval yields $|\Phi_4(s)|\lesssim s^{-3/2}$. Thus, with $\beta_2=1/2$ and $\beta_d=1$ for $d\ge3$, we conclude from \eqref{smult} that
\begin{equation}\label{bmult}
 |\widehat D_r(\xi)|_{\ell^2}\le C\min\left\{\frac{r|\xi|}{\sqrt d},\left(\frac{\sqrt d}{r|\xi|}\right)^{\beta_d}\right\}.
\end{equation}
Integrating \eqref{bmult} in the radius and applying Plancherel's theorem, we obtain
\begin{equation}\label{bshel}
 \int_{\mathbb R^d}\int_0^\infty |D_ra(x)|_{\ell^2}^2\,\frac{dr}{r}\,dx
 \le C\|a\|_2^2.
\end{equation}
The constant is uniform since $\int_0^\infty\min\{s,s^{-\beta_d}\}^2\,ds/s\le 3/2$.

To estimate a maximal tail using \eqref{bmult}, we localize in the spatial frequency. Let $\mathcal P_l$ be the Fourier projection onto $\sqrt d\,2^{-l-1}<|\xi|\le\sqrt d\,2^{-l}$ for $l\in\mathbb Z$. For every finite sequence $(u_l)$ with $u_l=\mathcal P_lu_l\in L^2$, we have
\begin{equation}\label{bpre}
 \left\|\sup_{j\in\mathbb Z}\left|\sum_{l\ge j}u_l\right|\right\|_2
 \le C\left(\sum_l\|u_l\|_2^2\right)^{1/2}.
\end{equation}
Compare the $j$th tail with $P_{2^j/\sqrt d}\sum_lu_l$. On the $l$th annulus the two multipliers differ by at most $C2^{-|j-l|}$. By orthogonality in $l$ and summation in $j$, the differences have $L^2(\ell^2_j)$ norm at most $C(\sum_l\|u_l\|_2^2)^{1/2}$. The Poisson maximal inequality \eqref{apois} controls the comparison term by the same quantity, proving \eqref{bpre}. For vector-valued functions we may apply \eqref{bpre} in each coordinate and sum the squares.

\begin{proposition}\label{bdens}
Suppose $d\ge2$ and $0\le h\in L^2(\mathbb R^d)$. Let $h_r(z)$ be a jointly measurable real or complex family such that, outside one spatial null set,
\[
 |h_{r_0}(z)|+\sum_{j=1}^n|h_{r_j}(z)-h_{r_{j-1}}(z)|\le h(z)
 \qquad(0<r_0<\cdots<r_n,\ n\ge0).
\]
In particular, this condition holds if $0\le h_r\le h$ and $h_r$ is nondecreasing in $r$. Under the stated variation bound,
\begin{equation}\label{bmest}
 \left\|\left(\sum_{i=1}^d\sup_{t>0}
          \left|\int_t^\infty (D_rh_r)_i\,\frac{dr}{r}\right|^2
                \right)^{1/2}\right\|_2\le C\|h\|_2
\end{equation}
with an absolute constant.
\end{proposition}

\begin{proof}
Restrict $r$ temporarily to a finite union of consecutive dyadic annuli and write $r=2^kv$ with $1\le v<2$. For each $v$, the condition $2^kv>t$ selects a tail in $k$. By Minkowski's inequality it suffices to estimate the maximal tails of $\sum_kD_{r_k}h_{r_k}$ uniformly in $v$, where $r_k=2^kv$.

Decompose each $D_{r_k}h_{r_k}$ by the projections $\mathcal P_{k+m}$, $m\in\mathbb Z$. We first show that
\begin{equation}\label{bdual}
 \sum_k\|\mathcal P_{k+m}D_{r_k}h_{r_k}\|_{L^2(\ell^2)}^2
 \le C2^{-2\beta_d |m|}\|h\|_2^2 .
\end{equation}
Pair the sequence $(\mathcal P_{k+m}D_{r_k}h_{r_k})_k$ with a finite sequence $(V_k)$ in $L^2(\ell^2)$. Writing $\psi_k=D_{r_k}^*\mathcal P_{k+m}V_k$, we obtain the pairing $\int\sum_k h_{r_k}(z)\overline{\psi_k(z)}\,dz$. Let $k_0$ be the first index, and set $a_{k_0}=h_{r_{k_0}}$ and $a_j=h_{r_j}-h_{r_{j-1}}$ for $j>k_0$. Summing by parts, we have
\begin{equation*}
 \sum_k h_{r_k}\overline{\psi_k}=\sum_j a_j\sum_{k\ge j}\overline{\psi_k}.
\end{equation*}
Since $\sum_j|a_j|\le h$, the pairing has modulus at most
$\|h\|_2\|\sup_j|\sum_{k\ge j}\psi_k|\|_2$.
The Fourier support of $\psi_k$ lies in the $(k+m)$th annulus. There the multiplier of $D_{r_k}^*$, viewed as a linear functional on $\ell^2$, has norm at most $C2^{-\beta_d|m|}$ by \eqref{bmult}. Applying \eqref{bpre} to $(\psi_k)$, we bound the pairing by
\begin{equation*}
 C2^{-\beta_d|m|}\|h\|_2
       \left(\sum_k\|V_k\|_{L^2(\ell^2)}^2\right)^{1/2}.
\end{equation*}
This proves \eqref{bdual} by duality.

We now apply \eqref{bpre} to the tails of $(\mathcal P_{k+m}D_{r_k}h_{r_k})_k$ in each coordinate. Summing the squared bounds over the coordinates and using \eqref{bdual}, we obtain an $L^2$ bound of $C2^{-\beta_d|m|}\|h\|_2$ for the coordinatewise maximal tails. The bounds are summable in $m$ uniformly in $d$, since $\beta_d\ge1/2$. Each frequency decomposition converges in $L^2$ by Plancherel, and the range of $k$ is finite. We may therefore sum in $m$ and integrate in $v$ against $dv/v$ to conclude \eqref{bmest} for the restricted radii.

To pass to all radii, note that for each $t>0$, polar coordinates and $|h_r|\le h$ bound $\int_t^\infty|(D_rh_r)_i(x)|\,dr/r$ by
\begin{equation*}
 c_d\int_{|x-z|>t}|x-z|^{-d}h(z)\,dz<\infty
 \qquad(t>0).
\end{equation*}
The last integral is finite by Cauchy--Schwarz. Consequently, for almost every $x$ the tails are absolutely convergent and continuous in $t>0$. Let the finite unions of annuli increase to $(0,\infty)$. Fatou's lemma for maxima over finitely many rational $t$, followed by exhaustion of the positive rationals, proves \eqref{bmest}.
\end{proof}

For the constant family $h_r=a$, we may take $h=|a|$. By polar coordinates, \eqref{bmest} then becomes
\begin{equation}\label{bLtwo}
 \|\mathcal B_d a\|_2\le C\|a\|_2\qquad(a\in L^2).
\end{equation}
For $d=1$ the same estimate is the classical $L^2$ inequality for the maximal Hilbert transform (see e.g. \cite[Corollary~3.8.2]{Gra24}). In dimensions $d\ge2$ we also need the derivative of the truncation with respect to the logarithmic radius. For $g\in L^1\cap L^2$, polar coordinates in \eqref{rdef} and the shell formula \eqref{shell} identify the derivative of this profile. By \eqref{bshel}, we then have
\begin{equation}\label{bder}
 A_i(u,x):=R_{i,e^u}g(x),\qquad \partial_u A(u,x)=-D_{e^u}g(x),\qquad
 \|\partial_u A\|_{L^2(dx\,du;\ell^2)}\le C\|g\|_2 .
\end{equation}
The derivative identity holds almost everywhere and in distributions.

\section{The contact set and conjugate Poisson integrals}\label{s4}
For bounded nonnegative data, we apply the fractional obstacle decomposition of Ouyang, Spector, and Stockdale \cite[Theorem~1.2]{OSS26}. Let $0\le f\in L^1\cap L^\infty$, $M:=\|f\|_1$, and $\Lambda>0$. We use the fractional obstacle at height $\Lambda$ to control the Poisson extension of $f\mathbf1_{\{U>0\}}$ on the contact set $\{U=0\}$. A weighted energy identity will then allow us to take a separate supremum in height for each conjugate Poisson integral.

\subsection{The contact trace}

\begin{proposition}\label{cobs}
In every dimension there is a nonnegative $U\in H^1\cap L^1\cap C_0$. Write $\Omega:=\{U>0\}$, $F:=\mathbb R^d\setminus\Omega$, and $g:=f\mathbf1_\Omega$. There exists $h\in L^1\cap L^\infty$ such that
\begin{equation}\label{cdec}
 |D|U=g-h,\qquad 0\le h\le\Lambda,\qquad h=\Lambda\quad\text{a.e. on }\Omega.
\end{equation}
Moreover, $\int h=\int g$, $|\Omega|\le M/\Lambda$, and $f\le\Lambda$ almost everywhere on $F$. At almost every $x\in F$,
\begin{equation}\label{cpoi}
 \sup_{s>0}P_sg(x)\le \Lambda+2(1+3/d)^{-(d+3)/2}h(x)
                 \le (1+2e^{-3/2})\Lambda .
\end{equation}
\end{proposition}

\begin{proof}
Since $f\in L^1\cap L^2$, by \cite[Theorem~1.2]{OSS26} the fractional obstacle decomposition at order one  supplies nonnegative functions $U\in H^1\cap L^1$ and $\mu\in L^1\cap L^\infty$ with $|D|U=f-\mu=:q$. Moreover, $0\le\mu\le\Lambda$, $\int\mu=M$, and $\mu=\Lambda$ almost everywhere on $\{U>0\}$. Integrating $\partial_vP_vU=-P_vq$ in $L^2$ yields, for $s>0$,
\[
 U=P_sU+\int_0^sP_vq\,dv .
\]
For $v>0$ we have $P_vq\in C_0$, and the contribution from $0<v<\epsilon$ has supremum norm at most $\epsilon\|q\|_\infty$. The right side thus represents $U$ by a function in $C_0$. Since the continuous representatives for different $s$ agree everywhere, this determines a single closed contact set $F$.

At a Lebesgue point $x\in F$ of $q$, we have $U(x)=0$ and $s^{-1}\int_0^sP_vq(x)\,dv\to q(x)$. Hence $P_sU(x)/s\to-q(x)$ as $s\downarrow0$. In the kernel for $P_sU(x)/s$ the integrand increases as $s$ decreases, so
\begin{equation}\label{ctr}
 0\le \frac{P_sU(x)}s\le c_d\int\frac{U(y)}{|x-y|^{d+1}}\,dy
     =-q(x)=\mu(x)-f(x).
\end{equation}
In particular $h:=\mu-f\mathbf1_F$ has all the properties in \eqref{cdec}, with $\int h=\int g$ and $\Lambda|\Omega|\le\int\mu$. At each such point we may also differentiate the Poisson kernel at positive height to obtain
\[
 P_s(g-h)(x)=-\partial_sP_sU(x)
  =c_d\int\frac{U(y)}{|x-y|^{d+1}}
       \frac{d\,a-1}{(1+a)^{(d+3)/2}}\,dy,
 \qquad a=\frac{s^2}{|x-y|^2}.
\]
The scalar factor attains its maximum $2(1+3/d)^{-(d+3)/2}$ at $a=3/d$. Since $P_sh\le\Lambda$, \eqref{ctr} implies \eqref{cpoi}.
\end{proof}

By the Fourier multipliers of $\mathbf R$ and $|D|$, the decomposition \eqref{cdec} implies
\begin{equation}\label{rcan}
 \mathbf R(g-h)=\mathbf R|D|U=-\nabla_x U\quad\text{in }L^2(\mathbb R^d;\ell^2).
\end{equation}
Since $\nabla_xU=0$ almost everywhere on $F=\{U=0\}$, the principal values of $\mathbf Rg$ and $\mathbf Rh$ agree there by \eqref{rcan}. To estimate their extensions at positive heights, we need a bound on $\nabla_xP_sU$.

\subsection{Weighted energy and conjugate Poisson integrals}

Suppose $Y$ satisfies
\begin{equation}\label{exteq}
 Y_{ss}+(1-m)s^{-1}Y_s+\Delta_xY=0 .
\end{equation}
For a smooth convex function $\Psi$, put $J(s)=\int\Psi(Y(s,x))\,dx$. Applying the chain rule to $\Psi(Y)$ and integrating by parts in $x$, we obtain
\begin{equation}\label{enid}
 (sJ'(s)-mJ(s))'=s\int\Psi''(Y)|\nabla_{x,s}Y|^2\,dx .
\end{equation}
In both applications, $Y$ and its first two derivatives admit integrable spatial majorants on every compact interval of positive heights. The derivatives $\Psi'$ and $\Psi''$ are bounded on $[0,\infty)$, so one may justify the integration by parts with spatial cutoffs. Suppose that $J$ has finite limits at zero and infinity. Applying the mean value theorem to $u\mapsto J(e^u)$ on intervals of length one approaching each endpoint, we find sequences along which $sJ'(s)\to0$. Integrating \eqref{enid} along these sequences, we see that the total weighted energy equals $m(J(0+)-J(\infty))$.

For $H>0$ we use the two convex functions
\[
 \Phi_H(v):=Hv+H^2(e^{-v/H}-1),\qquad \Psi_H(v):=\Phi_H(v)-Hv .
\]
They have the same second derivative $e^{-v/H}$, bounded below by $e^{-1}$ on $[0,H]$. For $m>0$ and $0\le a\in L^1\cap L^\infty$, direct differentiation of \eqref{aker} shows that $Y=W_s^{(m)}*a$ satisfies \eqref{exteq}. Take $J=\int\Phi_H(Y)$. As $\Phi_H$ is $H$-Lipschitz on $[0,\infty)$ and $Y(s)\to a$ in $L^1$, we have $J(0+)=\int\Phi_H(a)\le H\|a\|_1$. At infinity, $\|Y(s)\|_\infty\to0$, $\|Y(s)\|_1=\|a\|_1$, and $\Phi_H(v)=O(v^2)$ near zero. Hence $J(\infty)=0$. On $E_H:=\{x:\sup_s W_s^{(m)}*a(x)\le H\}$ this proves
\begin{equation}\label{bene}
 \int_{E_H}\int_0^\infty s|\nabla_{x,s}Y(s,x)|^2\,ds\,dx
 \le e\,mH\|a\|_1 .
\end{equation}

We denote the convolution kernel of $P_sR_i$ by $C_{s,i}$. In our Fourier convention,
\begin{equation}\label{cpker}
 P_sR_i=-\int_s^\infty\partial_{x_i}P_v\,dv,\qquad C_{s,i}(z)=c_d\frac{z_i}{(|z|^2+s^2)^{(d+1)/2}} .
\end{equation}
The operator identity follows by integrating the Fourier multiplier in $v$. Integrating the spatial derivative of the Poisson kernel in the same variable yields the expression for $C_{s,i}$.

For the operator $\mathcal N$ in \eqref{nmax}, we use the one-dimensional Sobolev inequality in the variable $u$. If $E\subset\mathbb R^d$ is measurable, then
\begin{equation}\label{logsp}
 \|\mathbf1_E\mathcal N V\|_2^2
 \le 2\|\mathbf1_E V\|_{L^2(dx\,du;\ell^2)}
      \|\mathbf1_E\partial_u V\|_{L^2(dx\,du;\ell^2)}
\end{equation}
whenever both factors are finite, with the continuous $H^1(du)$ representatives in each coordinate. Indeed, integration of $\partial_u|V_i|^2$ bounds $\sup_u|V_i|^2$ by $2\|V_i\|_{L^2(du)}\|\partial_uV_i\|_{L^2(du)}$. Summing over $i$ and applying Cauchy--Schwarz in $(i,x)$ proves \eqref{logsp}.

\begin{proposition}\label{csf}
For $g$ and $F$ as in Proposition~\ref{cobs},
\begin{equation}\label{csoft}
 \int_F\sum_{i=1}^d\sup_{s>0}|P_sR_i g(x)|^2\,dx\le C\Lambda M .
\end{equation}
\end{proposition}

\begin{proof}
To estimate $\nabla_xP_sU$ itself, take $V=P_sU/s$ and $J(s)=\int\Psi_\Lambda(V)$. Harmonicity of $P_sU$ implies that $V$ satisfies \eqref{exteq} with $m=-1$. As $s\downarrow0$, $V(s,x)$ increases to infinity on $\Omega$ and to $h(x)$ almost everywhere on $F$ by \eqref{ctr}. Since $0\le-\Psi_\Lambda(v)\le\min(\Lambda^2,\Lambda v)$,
\[
 -J(0+)=\Lambda^2|\Omega|+\Lambda^2\int_F(1-e^{-h/\Lambda})
 \le \Lambda\int h\le\Lambda M.
\]
At infinity $J(s)\to0$ because $\|V(s)\|_1=s^{-1}\|U\|_1$. On $F$ we have $\Psi_\Lambda''(V)\ge e^{-1}$ almost everywhere. The energy identity \eqref{enid} therefore yields
\begin{equation}\label{cene1}
 \int_F\int_0^\infty |\nabla_xP_sU|^2\,\frac{ds}{s}\,dx
 \le e\,[-J(0+)]\le e\Lambda M .
\end{equation}

For the derivative in $\log s$, we use $s\partial_s\nabla_xP_sU=-s\nabla_xP_s(g-h)$ from \eqref{cdec}. By \eqref{cpoi}, $P_*g\le2\Lambda$ almost everywhere on $F$, so we may apply \eqref{bene} to $g$ with $m=1$ and $H=2\Lambda$. For $h$, Plancherel's identity yields $\int_0^\infty s\|\nabla_xP_sh\|_2^2\,ds=\|h\|_2^2/4$. Consequently
\begin{equation}\label{cene2}
 \int_F\int_0^\infty|s\partial_s\nabla_xP_sU|^2\,\frac{ds}{s}\,dx\le C\Lambda M .
\end{equation}
With $u=\log s$, the two factors on the right of \eqref{logsp} for $\nabla_xP_{e^u}U$ are bounded by \eqref{cene1} and \eqref{cene2}, respectively. By \eqref{rcan} this controls the conjugate Poisson maxima for $g-h$. For the remaining density $h$, the Poisson maximal estimate \eqref{apois} and Plancherel's identity imply
\[
 \sum_i\big\|\sup_s|P_sR_i h|\big\|_2^2
 \le C\sum_i\|R_i h\|_2^2=C\|h\|_2^2\le C\Lambda M .
\]
This proves \eqref{csoft}.
\end{proof}

\subsection{Conjugate integrals and tangential gradients}
Now we derive analogous estimates for conjugate Poisson integrals and tangential derivatives of Poisson extensions, allowing the heights to vary independently. 

\begin{corollary}\label{cweak}
For every finite real or complex Borel measure $\nu$ on $\mathbb R^d$,
\[
 \left\|\left(\sum_{i=1}^d\sup_{s>0}|C_{s,i}*\nu|^2\right)^{1/2}
      \right\|_{L^{1,\infty}}\le C|\nu|(\mathbb R^d),
\]
where $C$ is absolute.
\end{corollary}

\begin{proof}
For a bounded nonnegative density $f\in L^1$, split $f=g+f\mathbf1_F$ at height $\Lambda=\lambda$. On $F$ we estimate the conjugate maxima for $g$ by \eqref{csoft}. For $f\mathbf1_F$ we may use \eqref{apois} and Plancherel's identity on all of $\mathbb R^d$ since $\|f\mathbf1_F\|_2^2\le\lambda M$. Chebyshev's inequality, together with $|F^c|\le M/\lambda$, proves the weak bound. Subadditivity extends it to bounded complex densities in $L^1$ by decomposition into the positive and negative parts of the real and imaginary parts.

For general $\nu$, let $\rho_\epsilon$ be a positive smooth mollifier of integral one. The functions $\rho_\epsilon*\nu$ are bounded and have $L^1$ norm at most $|\nu|(\mathbb R^d)$. For each $s>0$, the kernel $C_{s,i}$ is bounded and uniformly continuous, so $C_{s,i}*(\rho_\epsilon*\nu)\to C_{s,i}*\nu$ uniformly in $x$. We may thus pass to the limit for any finite set of positive rational heights. Exhausting those heights and using continuity in $s>0$ proves the assertion.
\end{proof}

\begin{corollary}\label{pgrad}
Let $d\ge1$. Suppose that $u\in L^1(\mathbb R^d)+L^2(\mathbb R^d)$ and that $|D|u=\nu$ is a finite complex Borel measure in the sense of distributions. Then
\begin{equation}\label{pgw11}
 \left\|\left(\sum_{i=1}^d\sup_{s>0}|\partial_{x_i}P_su|^2\right)^{1/2}
          \right\|_{L^{1,\infty}}\le C|\nu|(\mathbb R^d).
\end{equation}
One may take the constant of Theorem~\ref{main}.
\end{corollary}

\begin{proof}
Since $\widehat u\in L^\infty+L^2$, the identity $\widehat\nu=|\xi|\widehat u$ holds between locally integrable functions. For $s>0$, the Fourier transform of $C_{s,i}*\nu$ is $-i\xi_i e^{-s|\xi|}\widehat\nu(\xi)/|\xi|$. This function is integrable and equals the transform of $-\partial_{x_i}P_su$. We may therefore apply Corollary~\ref{cweak} to obtain \eqref{pgw11}.

To compare constants, define the radial factor $b_s$ by $C_{s,i}(z)=K_i(z)b_s(|z|)$ for $z\ne0$. As $b_s$ increases from zero to one, $b_s'$ is a probability density on $(0,\infty)$. Hence, by Fubini's theorem,
\begin{equation}\label{poim}
 b_s(r)=\left(\frac{r^2}{r^2+s^2}\right)^{(d+1)/2},\qquad
 C_{s,i}*\nu(x)=\int_0^\infty R_{i,t}\nu(x)b_s'(t)\,dt .
\end{equation}
An atom at $x$ contributes zero to both sides, and the integral of the absolute value is at most $\|C_{s,i}\|_\infty|\nu|(\mathbb R^d)$. By \eqref{poim}, the coordinatewise maxima are bounded in $\ell^2$ by $\mathcal B_d\nu$. The measure estimate in Corollary~\ref{mweak} proves the assertion about the constant.
\end{proof}

\section{Frequency estimates in the logarithmic radius}\label{s5}

We now return to the hard truncations in dimension $d\ge2$. For the source $g$ of Proposition~\ref{cobs}, let $A$ be the logarithmic truncation profile in \eqref{bder}. By \eqref{bLtwo}, $A(\cdot,x)$ is bounded for almost every $x$, while \eqref{bder} ensures that $\partial_uA\in L^2(du\,dx;\ell^2)$. Thus, for a bounded symbol $m_0$ with compact support away from zero, we define $m_0(D_u)A$ by its Fourier transform $m_0(\omega)\widehat{\partial_uA}/(i\omega)$.
For $V\in L^2(du\,dx;\ell^2)$, $b\in L^2(\mathbb R)$, and a measurable set $E\subset\mathbb R^d$,
Cauchy--Schwarz in the Fourier inversion formula, followed by Plancherel in $u$, yields
\begin{equation}\label{finv}
 \|\mathbf1_E\mathcal N[b(D_u)V]\|_2
 \le (2\pi)^{-1/2}\|b\|_2\|\mathbf1_E V\|_{L^2(du\,dx;\ell^2)} .
\end{equation}
In particular, the inverse Fourier integral is continuous in $u$ for almost every $x$.

Choose an even $p\in C_c^\infty(\mathbb R)$ supported in $[-2,2]$ and equal to one on $[-1,1]$, with $0\le p\le1$.
Let $F_0$ consist of the points of $F$ at which \eqref{cpoi} holds, and set $H=C_2L_d\Lambda$. Define
\begin{equation}\label{fset}
 F'=F_0\cap\{H_*g\le H\}
     \cap\bigcap_{1\le j,\ 2^j\le\sqrt d}
        \left\{x:\sup_{s>0}W_s^{(4^j)}*g(x)\le C_2(1+j)^3\Lambda\right\}.
\end{equation}
At $m=1$ we already have the Poisson bound \eqref{cpoi} on $F_0$. For $j\ge1$, the $j$th condition in \eqref{fset} removes a set of measure at most $CM/[\Lambda(1+j)^2]$ by the maximal estimate \eqref{amax} for $W_s^{(4^j)}$. Summing in $j$ and applying the heat maximal estimate \eqref{aheat}, we obtain
\begin{equation}\label{bexc}
 |(F')^c|\le CM/\Lambda,\qquad P_*g\le2\Lambda,\qquad H_*g\le H
 \quad\text{on }F'.
\end{equation}

For $0<w\le1$ the triangular kernel and its Fourier transform are
\begin{equation}\label{tri}
 k_w(u)=w^{-1}(1-|u|/w)_+,\qquad
 \widehat{k_w}(\omega)=\left(\frac{\sin(w\omega/2)}{w\omega/2}\right)^2.
\end{equation}
In particular $k_w\ge0$, $\int k_w=1$, and $|\widehat{k_w}(\omega)|\le C(1+w|\omega|)^{-2}$.

\begin{proposition}\label{bsmoo}
For $0<w\le1$, let $k_w$ be as in \eqref{tri} and $F'$ as in \eqref{fset}. Then
\begin{equation}\label{bbest}
 \|\mathbf1_{F'}\mathcal N(k_w*_u A)\|_2
 \le C\big(1+d^{-1/4}L_d^{3/2}+w^{-1}d^{-5/4}L_d^{1/2}\big)\sqrt{\Lambda M}.
\end{equation}
\end{proposition}

At frequencies near zero we use the conjugate Poisson estimate on $F$ from Section~\ref{s4}. At intermediate frequencies the choice of $m$ in $W_s^{(m)}$ depends on the band. For the high frequencies we pass to complex heat times and use the bound on $H_*g$ in \eqref{bexc}.

\subsection{Low and intermediate Mellin frequencies}

For $m\ge1$, let $\beta=(m+1)/2$ and
$b_m=\Gamma(m/2)/(\sqrt\pi\Gamma((m+1)/2))$. We integrate the spatial gradient of the corresponding extension in height and write
\begin{equation}\label{tgrad}
 \mathcal T_m(u)=-b_m\int_{\sqrt{\beta/\alpha}\,e^u}^\infty\nabla_x(W_s^{(m)}*g)\,ds.
\end{equation}
With $\rho=|x-z|$, we substitute $t=s^2/\rho^2$ in the differentiated kernel. The result is
$K(x-z)$ times the upper tail of the beta-prime density
$t^{\beta-1}(1+t)^{-\alpha-\beta}/\mathrm B(\beta,\alpha)$ starting at
$(\beta/\alpha)e^{2u}/\rho^2$. With $t=(\beta/\alpha)e^{2v}$, write
$\eta_m(v)=2t^\beta(1+t)^{-\alpha-\beta}/\mathrm B(\beta,\alpha)$. Thus
$\int\eta_m=1$ and $\mathcal T_m=\eta_m*_u A$. Indeed, the indicator inside $A(u-v,x)$ is $\mathbf1_{\{v>u-\log\rho\}}$. For fixed $u$, the upper tail is $O(\rho^{d+1})$ at the origin and at most one elsewhere. The kernel integration against $g\in L^1$ is therefore absolutely convergent. By the beta integral, the multiplier in the logarithmic radius is
\begin{equation}\label{betam}
 \widehat \eta_m(\omega)
 = (\beta/\alpha)^{i\omega/2}
   \frac{\Gamma(\beta-i\omega/2)\Gamma(\alpha+i\omega/2)}
        {\Gamma(\beta)\Gamma(\alpha)}.
\end{equation}
Taking logarithms in Euler's product \eqref{agam} and using $\log(1+t)\le t$ and $\sum_{n\ge0}(n+a)^{-2}\le2/a$ for $a\ge1$, we obtain $|\widehat\eta_m(\omega)|\ge\exp[-C\omega^2(1/m+1/d)]$.

Differentiating $\mathcal T_m$ in \eqref{tgrad} yields $\partial_u\mathcal T_m=b_m s\nabla_x(W_s^{(m)}*g)$ at $s=\sqrt{\beta/\alpha}\,e^u$. On the set $E_\kappa$ of \eqref{bene} with $a=g$ and $H=\kappa$, we may therefore apply the weighted energy bound. Since $m b_m^2\le C$ by the gamma ratio estimate, we find
\begin{equation}\label{smeng}
 \|\mathbf1_{E_\kappa}\partial_u \mathcal T_m\|_{L^2(du\,dx;\ell^2)}\le C\sqrt{\kappa M}.
\end{equation}
If $m_0$ is supported in $\{L\le|\omega|\le2L\}$, $L\ge1$, apply \eqref{finv} and \eqref{smeng} with
$b=m_0/(i\omega\widehat \eta_m)$ and $V=\partial_u \mathcal T_m=\eta_m*_u\partial_u A$. The lower bound for $|\widehat\eta_m|$ then implies
\begin{equation}\label{bmid}
 \|\mathbf1_{E_\kappa}\mathcal N[m_0(D_u)A]\|_2
 \le C L^{-1/2}e^{CL^2(1/m+1/d)}\|m_0\|_\infty\sqrt{\kappa M}.
\end{equation}
For $L=2^j\le\sqrt d$ take $m=L^2$. On $F'$ we may use $\kappa=C_2(1+j)^3\Lambda$ by \eqref{fset} and \eqref{cpoi}. The resulting bound in \eqref{bmid} is $C2^{-j/2}(1+j)^{3/2}\|m_0\|_\infty\sqrt{\Lambda M}$.

In the band containing zero we use \eqref{tgrad} at $m=1$. By \eqref{cpker}, $\mathcal T_{1,i}(u)=P_{e^u/\sqrt\alpha}R_ig$. Using the multiplier \eqref{betam} of $\eta_1$, write
\[
 p(D_u)A=p(D_u)\mathcal T_1+b(D_u)\partial_u \mathcal T_1,\qquad
 b(\omega)=\frac{p(\omega)(1-\widehat \eta_1(\omega))}
                   {i\omega\widehat \eta_1(\omega)}.
\]
Indeed $\int_v^\infty \eta_1=(1+e^{2v}/\alpha)^{-\alpha}$, which is at most
$(1+e^{2v})^{-1}$ for $v\ge0$. The lower tail for $v<0$ is at most $e^{2v}$.
These tail bounds imply a uniform bound on the first absolute moment of $\eta_1$ and finiteness of its higher moments. Consequently $b$ extends smoothly across zero, and the lower bound for $|\widehat\eta_1|$ on the support of $p$ yields $\|b\|_2\le C$. The decomposition holds in distributions because $p=p\widehat\eta_1+b\,i\omega\widehat\eta_1$ also at $\omega=0$. It thus accounts for any constant term in $A$.
Apply \eqref{csoft} and $\|\mathcal F^{-1}p\|_1<\infty$ to the first summand,
and use \eqref{finv} and the derivative bound \eqref{smeng} for $\partial_u\mathcal T_1$ with $\kappa=2\Lambda$ in the second. By \eqref{cpoi}, $F\subset E_{2\Lambda}$ up to a null set. We obtain
\begin{equation}\label{blow}
 \|\mathbf1_F\mathcal N[p(D_u)A]\|_2\le C\sqrt{\Lambda M}.
\end{equation}

\subsection{Complex heat estimates for high frequencies}

For $|\vartheta|<\pi/2$ let $c=\cos\vartheta$ and set
$Z(t)=c^{d/2}H_{te^{i\vartheta}}g$, $Y(t)=H_{t/c}g$.
The complex Gaussian satisfies $c^{d/2}|H_{te^{i\vartheta}}(x)|=H_{t/c}(x)$, so $|Z|\le Y$ for $g\ge0$. To bound the gradient of $Z$ on $F'$, put
\[
 b(y)=(1+y/H)^{-2},\qquad \phi(y)=Hy-H^2\log(1+y/H).
\]
Here $\phi''=b$, $b''\ge0$, $y|b'|\le2b$, $0\le\phi(y)\le Hy$ and $b(y)y^2\le Hy$.
For an absolute constant $C_0$ define
$J(t)=\int[C_0c^{-2}\phi(Y)+b(Y)|Z|^2]$. Then $0\le J(t)\le Cc^{-2}HM$.
Using $Y_t=c^{-1}\Delta Y$ and $Z_t=e^{i\vartheta}\Delta Z$, we differentiate $J$ and integrate by parts in $x$. The contribution of direction $j$ to $-J'(t)$ is the integral of
\[
 c^{-1}(C_0c^{-2}b+b''|Z|^2)q^2+
 2b'q\, Z\cdot (e^{i\vartheta}v+c^{-1}v)+2cb|v|^2 ,
\]
where $q=\partial_jY$, $v=\partial_j Z$, the dot is the real inner product on
$\mathbb C$, and $b,b',b''$ are evaluated at $Y$. The mixed term has modulus at most
$8c^{-1}b|q||v|\le cb|v|^2+16c^{-3}bq^2$. With $C_0=16$ we therefore have
$-J'\ge c\int b(Y)|\nabla Z|^2$.
On a compact positive time interval the Gaussian derivatives are integrable in
space and their convolutions with $g\in L^1\cap L^\infty$ are bounded and integrable.
Spatial cutoffs therefore justify the integration by parts. Since $b(Y)\ge1/4$ on $F'$,
we obtain upon integration over $0<t<\infty$
\begin{equation}\label{bheat}
 \int_{F'}\int_0^\infty|\nabla Z(t,x)|^2\,dt\,dx\le Cc^{-3}HM .
\end{equation}

Write $v=2u$, $a=e^v/(4\alpha)$ and $A^\natural(v)=A(v/2)$. Substitution
$r=|z|^2/(4t)$ shows that the heat tail
\begin{equation}\label{qheat}
 Q_{i,a}g=-\pi^{-1/2}\int_a^\infty\partial_{x_i}H_t g\,\frac{dt}{\sqrt t}
\end{equation}
has kernel $K_i(z)G_\alpha(|z|^2/(4a))$. Here
\begin{equation}\label{gdist}
 G_\alpha(b)=\int_0^b\gamma_\alpha(r)\,dr,\qquad
 \gamma_\alpha(r)=\frac{r^{\alpha-1}e^{-r}}{\Gamma(\alpha)}\quad (r>0).
\end{equation}
Thus $Q_a g=\ell_\alpha*_v A^\natural$ with $\ell_\alpha(v)=L_\alpha(\log\alpha-v)$ as in \eqref{subor}. Its Fourier transform is $\widehat\ell_\alpha(\tau)=\alpha^{-i\tau}\Gamma(\alpha+i\tau)/\Gamma(\alpha)$.
Put $\mathcal D_\vartheta(v)=\pi^{-1/2}e^{i\vartheta/2}\sqrt a\,\nabla H_{ae^{i\vartheta}}g$.
Then $\mathcal D_0=\partial_v Q_a g$. At each spatial frequency $\xi\ne0$, the gamma integral
\[
 e^{i\vartheta/2}\int_0^\infty
 a^{1/2-i\tau}e^{-ae^{i\vartheta}|\xi|^2}\,\frac{da}{a}
 =e^{-\vartheta\tau}\Gamma(1/2-i\tau)|\xi|^{-1+2i\tau}
\]
implies, after Fourier transformation in $v$,
\begin{equation}\label{mheat}
 \widehat{\mathcal D_\vartheta}(\tau,x)
 =e^{-\vartheta\tau}\widehat{\mathcal D_0}(\tau,x)
 =e^{-\vartheta\tau}\widehat{\ell_\alpha}(\tau)\,
       \widehat{\partial_v A^\natural}(\tau,x).
\end{equation}
For this calculation one may first assume that $\widehat g$ is smooth with compact support away from zero. For general $g\in L^2$, spatial Plancherel yields $\|\mathcal D_\vartheta\|_{L^2(dv\,dx;\ell^2)}^2=(2\pi c)^{-1}\|g\|_2^2$. Approximation extends \eqref{mheat} on each compact $\tau$ interval. By the bound \eqref{bheat} for $\nabla Z$, the squared norm of $c^{d/2}\mathcal D_\vartheta$
on $F'$ is at most $Cc^{-3}HM$. For a compact interval $I$ on one side of zero, let $\varepsilon$ be the sign of its midpoint and choose $\vartheta=-\varepsilon\theta$, $0\le\theta<\pi/2$, $c=\cos\theta$.
For a symbol $m_0$ supported in $I$, inversion of the multiplier in \eqref{mheat} and application of \eqref{finv} to $c^{d/2}\mathcal D_\vartheta$ lead to
\begin{equation}\label{bray}
 \|\mathbf1_{F'}\mathcal N[m_0(D_v)A^\natural]\|_2
 \le C\sqrt{HM}\|m_0\|_\infty c^{-3/2}
 \left(\int_I\left|\frac{c^{-d/2}e^{-\theta|\tau|}\Gamma(\alpha)}
  {\tau\Gamma(\alpha+i\tau)}\right|^2 d\tau\right)^{1/2}.
\end{equation}

\subsection{High-frequency bounds and summation}

It remains to choose the heat angle in \eqref{bray}. For $N=d/2$ and $\tau\ge0$, the minimum
\[
 T(\tau)=\min_{0\le\theta<\pi/2}(-N\log\cos\theta-\theta\tau)
 =\tfrac N2\log(1+\tau^2/N^2)-\tau\arctan(\tau/N)
\]
is attained at $\theta(\tau)=\arctan(\tau/N)$. We have
\begin{equation}\label{bgam}
 |\Gamma(\alpha+i\tau)|/\Gamma(\alpha)\ge e^{T(\tau)} .
\end{equation}
Indeed $x\mapsto\log(1+\tau^2/x^2)$ is convex. By \eqref{agam} and the midpoint
inequality on the successive unit intervals starting at $N$, the logarithm of the
left side is bounded below by $-\tfrac12\int_N^\infty\log(1+\tau^2/x^2)\,dx=T(\tau)$.
With $\theta=\theta(\tau_0)$, the affine function $\tau\mapsto-N\log\cos\theta-\theta\tau$ is tangent to $T$ at $\tau_0$. Since $T''(\tau)=-N/(N^2+\tau^2)$, this line lies at most $C$
above $T$ whenever $|\tau-\tau_0|\le c_0\sqrt{N^2+\tau_0^2}/\sqrt N$,
with $c_0$ a small absolute constant.

For a band $B\le|\omega|\le2B$ in $u$ with $B\ge\sqrt d$, the positive band in $v$
is $B/2\le\tau\le B$. Split it into $n$ equal intervals of length at most
$c_1\max(\sqrt d,B/\sqrt d)$, with $n$ minimal. If $\tau_I$ is the midpoint of $I$, take $\theta_I=\arctan(\tau_I/N)$. The tangent estimate and \eqref{bgam} bound the integrand before squaring in \eqref{bray} by $C/B$. The sum of the bounds over the intervals is at most
\[
 C\frac{\sqrt{nB}}{B}\max_I(\cos\theta_I)^{-3/2}
 \|m_0\|_\infty\sqrt{HM},
\]
since $\sum_I |I|^{1/2}\le(nB)^{1/2}$. Changing the signs of the angles yields the same bound on the negative band. For $\sqrt d\le B\le d$, we have $n\lesssim B/\sqrt d$ and $\cos\theta_I\gtrsim1$. For $B>d$, we have $n\lesssim\sqrt d$ and $\cos\theta_I\asymp d/B$. Since $\omega=2\tau$, the estimate in the original variable $u$ becomes
\begin{equation}\label{bhigh}
 \|\mathbf1_{F'}\mathcal N[m_0(D_u)A]\|_2
 \le C(d^{-1/4}+B d^{-5/4})\|m_0\|_\infty\sqrt{HM}.
\end{equation}

\begin{proof}[Proof of Proposition~\ref{bsmoo}]
Since $k_w$ in \eqref{tri} has integral one, the low-frequency estimate \eqref{blow} applies also to $k_w*_u p(D_u)A$.
For the multiplier $(1-p)\widehat{k_w}$, use \eqref{tri} and sum \eqref{bmid} with $m=B^2$ over the dyadic bands $B=2^j\le\sqrt d$. The coefficients $2^{-j/2}(1+j)^{3/2}$ are summable. For $\sqrt d<B\le d$, the estimate \eqref{bhigh} bounds each of the $O(L_d)$ bands by $Cd^{-1/4}\sqrt{HM}$. For $B>d$, use
$\sum_{B>d} B(1+wB)^{-2}\le C/w$, summing the two geometric series
on either side of $1/w$. Substitution of $H=C_2 L_d\Lambda$ proves \eqref{bbest}.
The partial sums outside the low band converge in $L^2(du\,dx;\ell^2)$
by Plancherel applied to $\partial_u A$. Their maximal norms are summable on $F'$,
so they also converge uniformly in $u$ for almost every $x\in F'$.
The two limits agree almost everywhere in $(u,x)$ on $\mathbb R\times F'$. Continuity in $u$ identifies the uniform limit there with $k_w*_u A$ for every $u$, for almost every $x\in F'$.
\end{proof}

\section{Variational equation and annular comparison}\label{s6}

To pass from the smoothed profile to the hard cutoff, we compare the Riesz kernels centered at a point $y$ of the source $g$ and at a point $z_*$ of the bounded density $h$. We construct a measure $\Pi$ with these marginals and estimate the conditional mass of $\{|y-z_*|\ge p/B\}$ in terms of the contact set on dyadic annuli about $y$. The comparison uses the even Poisson extension and the maximum principle on nested balls.

\subsection{The positive solution operator}

We work with the real Dirichlet form
\[
 \mathcal E(v,w)=\langle |D|^{1/2}v,|D|^{1/2}w\rangle
 =\frac{c_d}2\iint\frac{(v(x)-v(y))(w(x)-w(y))}{|x-y|^{d+1}}\,dx\,dy
 \qquad(v,w\in H^{1/2}).
\]
For $v=w$, the kernel representation follows by monotone convergence in $s^{-1}\langle(I-P_s)v,v\rangle$ as $s\downarrow0$. Polarization yields the general identity. Assume $\Omega\ne\varnothing$, since $g=h=0$ otherwise. Set $V:=\Lambda/U$ on $\Omega$, so $V\ge\Lambda/\|U\|_\infty>0$. The real space
\begin{equation*}
 \mathcal H:=\{v\in H^{1/2}: v=0\text{ a.e. on }F,\
                       \sqrt V v\in L^2(\Omega)\}
\end{equation*}
is a Hilbert space with squared norm $\mathcal E(v,v)+\int_\Omega V|v|^2$. Indeed a Cauchy sequence converges in $H^{1/2}$ by the lower bound on $V$, and its limit agrees almost everywhere on $\Omega$ with the limit in $L^2(V\,dx)$. In particular, $U\in\mathcal H$ since $\int_\Omega VU^2=\Lambda\int U$. Testing \eqref{cdec} against $v\in\mathcal H$ and using $VU=\Lambda=h$ on $\Omega$, we obtain
\begin{equation}\label{cvar}
 \mathcal E(U,v)+\int_\Omega VUv=\int g v\qquad(v\in\mathcal H).
\end{equation}

\begin{proposition}\label{cpair}
There is a finite positive measure
\begin{equation}\label{cdpi}
 d\Pi(y,z_*):=g(y)\,dy\,d\pi_y(z_*)
\end{equation}
such that, for every Borel set $E\subset\mathbb R^d$,
\begin{equation}\label{margs}
 \Pi(E\times\mathbb R^d)=\int_E g(y)\,dy,\qquad
 \Pi(\mathbb R^d\times E)=\int_E h(z_*)\,dz_*.
\end{equation}
The measures $\pi_y$ have mass one for $g(y)\,dy$ almost every $y$. For each $\varphi\in C_c^\infty(\mathbb R^d)$ with $0\le\varphi\le1$, there is a function $b_\varphi\in H^{1/2}$ such that
\begin{equation}\label{ctrb}
 b_\varphi(y)=\int\varphi\,d\pi_y\quad\text{for }g(y)\,dy\text{ a.e. }y,
 \qquad b_\varphi=\varphi\quad\text{a.e. on }F.
\end{equation}
Its even Poisson extension $\widetilde b_\varphi(x,t):=P_{|t|}b_\varphi(x)$ belongs to $H^1_{\rm loc}(\mathbb R^{d+1})$, lies in $[0,1]$ almost everywhere, and is weakly subharmonic in every ball centered on $\mathbb R^d\times\{0\}$ whose equatorial disk is contained in $\{\varphi=0\}$.
\end{proposition}

\begin{proof}
Let $\varphi\in C_c^\infty$ be real. The functional $v\mapsto-\mathcal E(\varphi,v)$ is bounded on $\mathcal H$. Its Riesz representative $b_\varphi-\varphi\in\mathcal H$ satisfies
\begin{equation}\label{cform}
 \mathcal E(b_\varphi,v)+\int_\Omega V(b_\varphi-\varphi)v=0
                       \qquad(v\in\mathcal H).
\end{equation}
If $0\le\varphi\le1$, then $b_\varphi^-$ and $(b_\varphi-1)_+$ belong to $\mathcal H$. They vanish on $F$ and are bounded by $|b_\varphi-\varphi|$ on $\Omega$. By the difference representation of $\mathcal E$, the left side of \eqref{cform} is at most $-\|b_\varphi^-\|_{\mathcal H}^2$ for the first test and at least $\|(b_\varphi-1)_+\|_{\mathcal H}^2$ for the second. Hence $0\le b_\varphi\le1$. For arbitrary real $\varphi$, the tests $(\pm b_\varphi-\|\varphi\|_\infty)_+$ likewise show that $\|b_\varphi\|_\infty\le\|\varphi\|_\infty$.

Choose a countable rational vector space of smooth compactly supported tests dense in $C_0(\mathbb R^d)$, containing a dense subset of its nonnegative cone. At each common Lebesgue point $y$ of the corresponding $b_\varphi$, take the precise representatives. Then $\varphi\mapsto b_\varphi(y)$ is rational-linear, is bounded by $\|\varphi\|_\infty$ in absolute value, and is nonnegative on the chosen nonnegative tests. By density and the Riesz representation theorem (see e.g. Folland \cite[Theorem~7.17]{Fol99}), it extends to integration against a positive measure $\pi_y$ of mass at most one. Set $\pi_y=0$ on the exceptional null set. By uniform approximation, the representation in \eqref{ctrb} holds almost everywhere for each smooth test. For $\varphi\in C_0$, the function $y\mapsto\int\varphi\,d\pi_y$ is measurable. Approximation on open sets and the monotone class theorem therefore imply that $y\mapsto\pi_y(E)$ is measurable for every Borel $E$.

Test the variational equation \eqref{cform} for $b_\varphi$ with $U$ and use $b_\varphi-\varphi$ in the identity \eqref{cvar}. By symmetry,
\begin{equation*}
 \int g(b_\varphi-\varphi)=-\mathcal E(\varphi,U)=-\int(g-h)\varphi.
\end{equation*}
It follows that $\int g b_\varphi=\int h\varphi$ for every $\varphi\in C_c^\infty$. By uniqueness of finite measures, the second marginal of $\Pi$ is $h\,dz_*$. Taking compact cutoffs increasing to one, we obtain $\int g(y)\pi_y(\mathbb R^d)\,dy=\int h=\int g$. Since $\pi_y(\mathbb R^d)\le1$, we infer that $\pi_y(\mathbb R^d)=1$ for $g(y)\,dy$ almost every $y$. Integration of \eqref{cdpi} over $z_*$ proves the first identity in \eqref{margs}.

For subharmonicity take $b=b_\varphi$ and a nonnegative smooth test $\psi$ supported in a disk where $\varphi=0$. For $\epsilon>0$, $w:=\min(b,\epsilon\psi)$ belongs to $\mathcal H$, since it vanishes almost everywhere on $F$ and $0\le w\le |b-\varphi|$ on its support. With $w$ as test function in \eqref{cform}, we have $\mathcal E(b,w)=-\int_\Omega Vbw\le0$. Set $z:=\epsilon\psi-b$. Since $w=\epsilon\psi-z^+$ and $\mathcal E(z,z^+)\ge\mathcal E(z^+,z^+)$ by the difference kernel,
\[
 \epsilon\mathcal E(b,\psi)\le \mathcal E(b,z^+)
 \le\epsilon\mathcal E(\psi,z^+)-\mathcal E(z^+,z^+)
 \le\tfrac14\epsilon^2\mathcal E(\psi,\psi).
\]
Letting $\epsilon\downarrow0$ proves $\mathcal E(b,\psi)\le0$.
By Plancherel, the extension $\widetilde b$ with multiplier $e^{-|t||\xi|}$ belongs to $H^1$ on each strip of finite height and has trace $b$. Via \cite[Section~1, (1.1)--(1.2)]{CaSi07}, the Dirichlet-to-Neumann identity for the Poisson extension gives $\partial_t\widetilde b(0\pm)=\mp|D|b$ in $H^{-1/2}$. The jump therefore yields $-\Delta\widetilde b=2|D|b\,\delta_{\{t=0\}}$. For nonnegative tests in the specified balls, $\mathcal E(b,\psi)\le0$ implies $\Delta\widetilde b\ge0$ in the distributional sense. Moreover, $0\le\widetilde b\le1$ by positivity and unit mass of the Poisson kernel.
\end{proof}

\subsection{A sum over contact annuli}

For $d\ge2$ put $I_1(x):=c_d|x|^{1-d}/(d-1)$ and, for each $y\in\mathbb R^d$ and dyadic $p$, set $E_p(y):=F\cap\{x:p\le|x-y|\le2p\}$. For solutions of the fractional Laplace equation, Wiener estimates express decay at regular boundary points in terms of the capacity of the complement (see Bj\"orn \cite[Proposition~3.3]{Bj24}). For the interior centers $y\in\Omega$ considered here, the potentials $I_1*\nu_p$ of densities on $E_p(y)$ quantify the decrease of the tails of $\pi_y$. We need the dependence on $d$ explicitly.

\begin{proposition}\label{cspr}
For $g(y)\,dy$ almost every $y$, consider any family of nonnegative $L^1$ densities $\nu_p$, indexed by the dyadic radii $p$, with $\nu_p=0$ almost everywhere off $E_p(y)$. Suppose that $I_1*\nu_p\le A$ almost everywhere with respect to $\nu_p(x)\,dx$, with the same $A>0$ for every $p$. Then, for every $B\ge1$,
\begin{equation}\label{csurv}
 \sum_{p} (I_1*\nu_p)(y)\,\pi_y\bigl\{z_*:|y-z_*|\ge p/B\bigr\}
 \le C A\bigl(L_d+\log(B+1)\bigr).
\end{equation}
\end{proposition}

\begin{proof}
Truncating each $\nu_p$ from above preserves the assumed potential bound, so it suffices by monotone convergence to treat bounded densities. We shall choose the exceptional set of centers independently of the densities. For now take $y\in\mathbb R^d$, write $Y=(y,0)$, and let $B_r(Y)$ denote the ball of radius $r$ about $Y$ in $\mathbb R^{d+1}$. For each $p$ define
\begin{equation*}
 W_p(x,t):=\frac{c_d}{(d-1)A}
          \int \frac{\nu_p(z)}{(|x-z|^2+t^2)^{(d-1)/2}}\,dz,
 \qquad v_p:=W_p(Y).
\end{equation*}
The kernel representation shows that $W_p$ tends to zero at infinity and is continuous. Indeed, the contribution from $|x-z|<\epsilon$ is at most $C_d A^{-1}\|\nu_p\|_\infty\epsilon$, uniformly in $t$, and away from this region continuity follows by dominated convergence. Also $W_p=A^{-1}\int_{|t|}^\infty P_s\nu_p\,ds$, whence $-\Delta W_p=2A^{-1}\nu_p\delta_{\{t=0\}}$. By Plancherel, its full gradient energy is at most $C A^{-2}\int |\widehat\nu_p(\xi)|^2|\xi|^{-1}\,d\xi$, which is finite for $d\ge2$. Thus $W_p\in H^1_{\rm loc}$. The trace bound and continuity imply $W_p\le1$ on the closed support of $\nu_p\,dx\,\delta_{\{t=0\}}$. Outside that compact set $W_p$ is harmonic and vanishes at infinity. The maximum principle now shows that $0\le W_p\le1$ everywhere.

Use the barrier $\gamma_p:=W_p-3^{1-d}v_p$. On $\partial B_{8p}(Y)$, every distance to $(z,0)$ in the integral is at least $3|y-z|$, so $\gamma_p\le0$ there. On $B_{p/(16d)}(Y)$ each distance is at most $(1+1/(16d))|y-z|$. It follows that $\gamma_p\ge c v_p$ throughout this smaller ball, with an absolute $0<c<1$. We also have $\gamma_p\le1$ and $-\Delta\gamma_p=2A^{-1}\nu_p\delta_{\{t=0\}}$.

Let $b:=\widetilde b_\varphi$ for a test as in Proposition~\ref{cpair} with $\varphi=0$ on the equatorial disk of $B_{8p}(Y)$, and put $m_p:=M_b(8p)$, where $M_b(R):=\operatorname*{ess\,sup}_{B_R(Y)} b$. The outer trace of $b$ is at most $m_p$, and $\gamma_p\le0$ there. Thus $\zeta:=(b-m_p(1-\gamma_p))_+$ belongs to $H^1_0(B_{8p}(Y))$ in view of \cite[Section~5.5, Theorem~2]{Eva10}. Nonnegative smooth approximation in $H^1_0$ permits $\zeta$ as a test for subharmonicity. Since $\nu_p$ is bounded and compactly supported, integration against $\nu_p\,dx$ on the equatorial disk is continuous on $H^1_0$ by the trace theorem. On that disk, $\varphi=0$ and $\gamma_p\le1$, so the trace relation \eqref{ctrb} implies that $\zeta=0$ at $\nu_p(x)\,dx$ almost every point. Testing the inequality for $b$ and the equation for $\gamma_p$ with $\zeta$, we obtain
\begin{equation*}
 \int_{B_{8p}(Y)}|\nabla\zeta|^2
 =\int_{B_{8p}(Y)}(\nabla b+m_p\nabla\gamma_p)\cdot\nabla\zeta\le0.
\end{equation*}
It follows that $b\le m_p(1-\gamma_p)$ almost everywhere in the ball and hence
\begin{equation}\label{ccon}
 M_b\bigl(p/(16d)\bigr)\le(1-cv_p)M_b(8p).
\end{equation}

Partition the dyadic radii $q=2^k$ according to the residue of $k$ modulo $J=1+\lceil\log_2(128d)\rceil$. Within each class the intervals $[q/(16d),8q]$ are disjoint and ordered. For a fixed $p$ and a test vanishing on the equatorial disk of $B_{8p}(Y)$, iterate the contraction \eqref{ccon} for $M_b$ over finite ranges of $q\le p$ in the class of $p$, using monotonicity of $M_b$. At a Lebesgue point $y$ of $b_\varphi$, $P_t b_\varphi(y)\to b_\varphi(y)$ as $t\downarrow0$. At positive heights the extension is continuous, so the bound by the essential supremum on the innermost ball passes to this limit. Letting the lower scale tend to zero gives
\begin{equation*}
 b_\varphi(y)\le \prod_{\substack{q\le p\\q\text{ in the class of }p}}(1-cv_q).
\end{equation*}
Here we used $M_b(8p)\le1$.

To apply the product bound to $\pi_y$ outside a single null set of centers, it suffices to use countably many tests. Choose a smooth function $0\le\eta\le1$ equal to zero on $[0,10]$ and to one on $[14,\infty)$, and a smooth compactly supported cutoff $0\le\chi\le1$ on $\mathbb R^d$ equal to one near zero. Use the tests $\eta(|z-a|/p)\chi(z/j)$ with $a\in\mathbb Q^d$, $p$ dyadic and $j$ a positive integer. On a common set of full $g(y)\,dy$ measure, $y$ is a Lebesgue point for each trace and the representation \eqref{ctrb} holds for all these tests. For $y$ in that set choose $|a-y|<p$. Each test vanishes on $|z-y|<8p$, and its limit as $j\to\infty$ equals one on $|z-y|\ge16p$. Thus
\begin{equation*}
 \pi_y\{|z_*-y|\ge16p\}\le \exp\left(-c\sum_{\substack{q\le p\\q\text{ in the class of }p}}v_q\right).
\end{equation*}
The exceptional set depends only on the tests. On a finite range of radii in one class, let $T_p$ be the sum of the $v_q$ over that range up to $p$. Since $v_p e^{-cT_p}\le\int_{T_p-v_p}^{T_p}e^{-ct}\,dt$, summation in the class costs at most $1/c$, independently of the range. Exhausting each class and summing over the $J$ classes, we obtain $\sum_p v_p\pi_y\{|z_*-y|\ge16p\}\le CL_d$. For each $z_*$, at most $C\log(B+1)$ dyadic radii satisfy $p/B\le|z_*-y|<16p$. Since $v_p\le1$ and $\pi_y$ has mass one, the contribution of these radii is bounded by $C\log(B+1)$ after integration in $z_*$. Substituting $v_p=(I_1*\nu_p)(y)/A$ proves \eqref{csurv}.
\end{proof}

\section{Passage from smoothed to hard truncations}\label{s7}

For $0\le f\in L^1\cap L^\infty$, let $g,h$ be the densities from Proposition~\ref{cobs} and let $\Pi$ be their coupling from Proposition~\ref{cpair}. For an observation point $x$ and a pair $(y,z_*)$, write
\[
 s=|x-y|,\qquad R=|x-z_*|,\qquad P=|y-z_*|,\qquad
 \delta=d^{-2}.
\]
If $g=0$, then $f\le\Lambda$ almost everywhere and $\|f\|_2^2\le\Lambda M$, so \eqref{bLtwo} suffices. We may therefore assume that $g\ne0$. We first work in sufficiently large dimension.

\subsection{Change of kernel center}

When $P\le\delta R$, the distance between the two sources is small compared with their distances to $x$. We write the contribution of these pairs as
\[
 A_i^{\mathrm{far}}(u,x)
 =\int_{P\le\delta R}K_i(x-y)\mathbf1_{\{s>e^u\}}\,d\Pi(y,z_*).
\]
Define $A^{\mathrm{near}}$ by integrating over the complementary region. Since the first marginal of $\Pi$ is $g(y)\,dy$ by \eqref{margs}, we have $A=A^{\mathrm{far}}+A^{\mathrm{near}}$. We shall construct nonnegative functions $E_1,E_2$ such that
\begin{equation}\label{pfar}
 \mathcal N A^{\mathrm{far}}\le E_1+E_2,\qquad
 \|E_1\|_1\le CM,\qquad \|E_2\|_2^2\le C\Lambda M .
\end{equation}

To move the kernel center from $y$ to $z_*$, observe that every point on the segment joining $x-y$ to $x-z_*$ has norm at least $(1-\delta)R$. Since $(1-\delta)^{-d-1}\le C$, differentiation of the vector kernel yields
\begin{equation*}
 |K(x-y)-K(x-z_*)|_{\ell^2}\le Cd(P/R)c_dR^{-d}.
\end{equation*}
For $P>0$, integration in $x$ over $R\ge P/\delta$ bounds the right-hand side by $Cd\alpha_d\delta\le C$. The difference vanishes when $P=0$. Taking absolute values in each coordinate and then applying Minkowski's inequality, we obtain a majorant $E_1$ of the maximal error with $\|E_1\|_1\le CM$.

We next move the center of the cutoff ball. For this purpose we disintegrate $\Pi$ over its second marginal $h(z_*)\,dz_*$. Since both marginals are measures on $\mathbb R^d$, Kallenberg's disintegration theorem applies to the probability measure $\Pi/\int g$ by \cite[Theorem~8.5]{Kal21}. We may thus write
\begin{equation}\label{hmarg}
 d\Pi(y,z_*)=d\kappa_{z_*}(y)\,h(z_*)\,dz_*,
 \qquad h_r(z_*):=h(z_*)\kappa_{z_*}\{y:|y-z_*|\le r\},
\end{equation}
where $\kappa$ is a Borel probability kernel. Thus $h_r(z_*)$ is jointly measurable in $(r,z_*)$, increases with $r$, and is bounded by $h(z_*)$. In particular, we can use \eqref{hmarg} at a radius depending measurably on $(x,z_*)$.

Write $t=e^u$ and $q_j=2^{-j}\delta$ for $j\ge0$. On $q_j/2<P/R\le q_j$, let $G_{i,j}(t,x)$ denote the error incurred by replacing $s>t$ with $R>t$ in the integral with kernel $K_i(x-z_*)$. The two cutoffs agree when $P=0$. On the $j$th range they can differ only if $|t-R|\le q_jR$, since $|s-R|\le P$. Put
\begin{equation*}
 h_{R,j}:=h_{q_jR}-h_{q_jR/2},\qquad
 S_{R,i}a(x):=\int_{S^{d-1}}|\theta_i|a(x-R\theta)\,d\sigma_{d-1}(\theta).
\end{equation*}
In polar coordinates centered at $x$, the absolute cutoff error therefore satisfies
\begin{equation*}
 |G_{i,j}(t,x)|
 \le \alpha_d\int_{t/(1+q_j)}^{t/(1-q_j)}
          S_{R,i}h_{R,j}(x)\,\frac{dR}{R}.
\end{equation*}
The interval has $dR/R$ measure $\log((1+q_j)/(1-q_j))\le Cq_j$. After Cauchy--Schwarz in $R$, we may extend the integral of the square to all $R>0$ and take the supremum in $t$ separately for each coordinate. Using $|(S_{R,i}a)_i|_{\ell^2}\le S_Ra$ for $a\ge0$, we obtain
\[
 \sum_i\sup_{t>0}|G_{i,j}(t,x)|^2
 \le C\alpha_d^2q_j\int_0^\infty(S_Rh_{R,j}(x))^2\,\frac{dR}{R}.
\]
These estimates hold simultaneously for rational $t>0$ outside a single null set. Absolute continuity of both marginals ensures continuity of the truncated integrals in $t>0$.

To integrate in $x$, note that for each $P>0$ the set of $R$ such that $q_j/2<P/R\le q_j$ has $dR/R$ measure $\log 2$. The set is empty for $P=0$. Hence, by \eqref{hmarg} and Tonelli's theorem, $\int_0^\infty h_{R,j}\,dR/R\le(\log 2)h$ almost everywhere. Since $0\le h_{R,j}\le h$ and $S_R$ is a contraction on $L^2$,
\begin{equation*}
 \int_0^\infty\|S_Rh_{R,j}\|_2^2\,\frac{dR}{R}
 \le \int_{\mathbb R^d}h(z_*)\int_0^\infty h_{R,j}(z_*)\,\frac{dR}{R}\,dz_*
 \le C\|h\|_2^2.
\end{equation*}
Since $\alpha_d\sum_j\sqrt{q_j}\lesssim\sqrt{d\delta}\le1$, the bounds are summable in $j$. Thus the $\ell^2$ norm of the coordinatewise maximal cutoff errors has $L^2$ norm at most $C\|h\|_2$.

With both centers at $z_*$, the shell formula \eqref{shell} reads
\begin{equation*}
 \int_{\{P\le\delta R\}}K_i(x-z_*)\mathbf1_{\{R>e^u\}}\,d\Pi(y,z_*)
   =\int_{e^u}^\infty (D_R h_{\delta R})_i(x)\,\frac{dR}{R}.
\end{equation*}
Since $h_{\delta R}$ increases with $R$ and is bounded by $h$, the same $L^2$ bound for the coordinatewise maxima follows from Proposition~\ref{bdens}. Together with $\|h\|_2^2\le\Lambda M$, this proves \eqref{pfar}.
On the complementary region $P>\delta R$, let $\theta=(x-y)/s$ and write
\begin{equation}\label{near}
 d\mu_x(y,z_*)=c_ds^{-d}\mathbf1_{\{R<d^2P\}}\,d\Pi(y,z_*),
 \qquad A_i^{\mathrm{near}}(u,x)=\int\theta_i\mathbf1_{\{\log s>u\}}\,d\mu_x.
\end{equation}

The bound $\mu_x\{s>a\}\le c_da^{-d}\|g\|_1$ for $a>0$ ensures that the tails in \eqref{near} can be averaged over bounded intervals of logarithmic radii.

\subsection{Cancellation between two sources}

Let $0<w\le c/d$ for a sufficiently small absolute constant $c$. We compare the hard cutoff at $u$ with two translates of its average against the triangular kernel $k_w$ in \eqref{tri}. On $F$ put
\begin{equation}\label{perr}
 \begin{aligned}
 A^c(u,x)&=(k_w*_u A^{\mathrm{near}})(u,x),\\
 E_w(x)&=\left(\sum_i
 \left[\sup_u\left(|A_i^{\mathrm{near}}(u,x)|
              -\max_{\pm}|A_i^c(u\pm w,x)|\right)_+\right]^2
 \right)^{1/2}.
 \end{aligned}
\end{equation}
By construction, $\mathcal N A^{\mathrm{near}}\le\mathcal N A^c+E_w$. Moreover, averaging against the probability density $k_w$ does not increase $\mathcal N A^{\mathrm{far}}$. Writing $A^c=k_w*_u A-k_w*_u A^{\mathrm{far}}$ and applying \eqref{pfar}, we find
\begin{equation}\label{recov}
 \mathcal B_dg=\mathcal N A
 \le \mathcal N(k_w*_u A)+2(E_1+E_2)+E_w\qquad\text{a.e. on }F.
\end{equation}
It remains to bound the excess $E_w$ on $F'$. We claim that
\begin{equation}\label{pthin}
 \|\mathbf1_{F'}E_w\|_2^2
  \le Cdw L_d^2\big[1+\log(1/w)+L_d\big]\Lambda M.
\end{equation}

For each $i,x,u$, write $C_0:=A_i^c(u+w,x)$,
$B:=A_i^{\mathrm{near}}(u,x)-C_0$, and
$D_0:=A_i^c(u-w,x)-A_i^{\mathrm{near}}(u,x)$.
By positivity and the support of $k_w$, each difference is an integral of $\theta_i$ against a weight between zero and one with respect to $\mu_x$. For $B$ this weight vanishes outside $u<\log s\le u+2w$, and for $D_0$ it vanishes outside $u-2w\le\log s\le u$. Since $g\ge0$, all three quantities are real and satisfy
\begin{equation*}
 \left(|C_0+B|-\max\{|C_0|,|C_0+B+D_0|\}\right)_+
       \le \sqrt{(BD_0)_-}.
\end{equation*}
Indeed, when $BD_0\ge0$, the value $C_0+B$ lies between $C_0$ and $C_0+B+D_0$. When $BD_0<0$, the triangle inequality bounds the excess by $\min(|B|,|D_0|)$.

To estimate the negative product, let $z$ denote the source in the second factor, and set $s':=|x-z|$ and $\theta':=(x-z)/s'$. Then
\begin{equation*}
 \sum_{i=1}^d (\theta_i\theta_i')_-
 \le \min\left\{1,\frac{|\theta-\theta'|^2}{4}\right\}
 \le \min\left\{1,\frac{|y-z|^2}{4ss'}\right\}.
\end{equation*}
Here we used Cauchy--Schwarz, $(ab)_-\le(a-b)^2/4$, and
$ss'|\theta-\theta'|^2=|y-z|^2-(s-s')^2$.
For each coordinate the two restrictions in $BD_0$ imply $|\log(s/s')|\le4w$, a condition independent of $u$. We may therefore bound the negative part by integrating the negative part of the product, discard the other restrictions on the second pair, and use its first marginal $g(z)\,dz$. The resulting estimate is uniform in $u$ for every $i$. Summing over $i$, we obtain
 \begin{align}
 E_w(x)^2&\le \iint g(z)c_d(s')^{-d}\mathbf1_{\{|\log(s/s')|\le4w\}}\notag\\
 &\quad\times\min\{1,|y-z|^2/(4ss')\}\,dz\,d\mu_x(y,z_*) .\label{ppair}
 \end{align}

We shall apply the annular estimate of Proposition~\ref{cspr} to all pairs in each dyadic shell. Write $D=|y-z|$ and introduce the density
\begin{align}
 \nu_{y,p}(x):={}&\frac{\mathbf1_{F'\cap\{p\le s<2p\}}(x)}p\int
 g(z)c_d(s')^{-d}\mathbf1_{\{|\log(s/s')|\le4w\}}\,
 \min\{1,D^2/(4ss')\}\,dz .\label{nudef}
\end{align}
The indicator in \eqref{nudef} restricts $x$ to $E_p(y)$ and $s'$ to $[p/2,3p]$, so $\nu_{y,p}$ is bounded for each $y,p$. We distinguish the two ranges $D\le K_0\sqrt w\,p$ and $D>K_0\sqrt w\,p$, taking $K_0$ to be a sufficiently large absolute constant. In the first range the angular factor is at most $Cw$. Moreover, the Poisson kernel satisfies
\begin{equation}\label{ppois}
 c_d u_0^{-d}\le C\sqrt d\,P_{p/\sqrt d}(u_0)
                  \qquad (p/2\le u_0\le3p).
\end{equation}
Indeed, $c_d u_0^{-d}/P_{p/\sqrt d}(u_0)=(u_0\sqrt d/p)(1+p^2/(du_0^2))^{(d+1)/2}$. Combining \eqref{ppois} with $P_*g\le2\Lambda$ on $F'$ from \eqref{bexc}, we see that the contribution of these pairs to \eqref{nudef} is at most $C\sqrt d\,w\Lambda\,\mathbf1_{E_p(y)}/p$. The radial monotonicity of $I_1$ and the layer cake representation imply that
\[
 \|I_1*\mathbf1_{E_p(y)}\|_\infty
 \le\int_{|z|\le2p}I_1(z)\,dz=\frac{2p\alpha_d}{d-1}.
\]
Thus the potential of this part of $\nu_{y,p}$ is bounded everywhere by $Cw\Lambda$, since $\alpha_d\lesssim\sqrt d$.

In the second range, put $\ell=D/p$. The triangle inequality $D\le s+s'$ yields $K_0\sqrt w<\ell\le5$. After translating and rotating, write
$y=(0,-D/2)$, $z=(0,D/2)$, and $x=(v,t)$, where
$v\in\mathbb R^{d-1}$. Since $2Dt=s^2-(s')^2$ and
$|s-s'|\le Cwp$, the domain of integration is contained in
\begin{equation*}
 |t|\le U_0,\qquad U_0:=Cwp/\ell .
\end{equation*}
By our choice of $K_0$, we have $U_0\le D/4$. Throughout this slab we use the bound $C\ell^2$ for the angular factor on the support of \eqref{nudef}.

To estimate the potential at a point $x_0=(v_0,\tau)$ with $p\le|x_0-y|\le2p$, write
\begin{equation*}
 r:=|x_0-y|,\quad r':=|x_0-z|,\quad
 a:=|\tau-t|,\quad b_0:=D/2-t,\quad
 R_0:=\big(|v_0|^2+(a+b_0)^2\big)^{1/2}.
\end{equation*}
For $|t|\le U_0$ we have
\begin{equation}\label{pgeom}
 \begin{gathered}
 b_0\asymp \ell p,\qquad |\tau|=\frac{|r^2-(r')^2|}{2D},\\
 R_0\ge r',\quad R_0\ge r-2|t|,\quad R_0\le Cp .
 \end{gathered}
\end{equation}
For the two lower bounds, observe that $a+b_0=\max\{|\tau-D/2|,|\tau-(2t-D/2)|\}$. The point $(0,2t-D/2)$, obtained by reflecting $z$ across the slice at height $t$, is at distance $2|t|$ from $y$.

On each such slice both kernels admit representations in terms of the Poisson kernel $P_h^{(d-1)}$ on $\mathbb R^{d-1}$. Namely,
\begin{align*}
 c_d|x-z|^{-d}&=\frac{c_d}{c_{d-1}b_0}P_{b_0}^{(d-1)}(v),\\
 I_1(x_0-x)&=\frac1\pi\int_a^\infty
       \frac{P_h^{(d-1)}(v_0-v)}{\sqrt{h^2-a^2}}\,dh .
\end{align*}
For the second identity, substitute $h^2=a^2+u^2$ and use
$c_d/(d-1)=\Gamma((d-1)/2)/(2\pi^{(d+1)/2})$.
By the Poisson semigroup identity in $\mathbb R^{d-1}$,
\begin{align}
 Q(t)&:=\int_{\mathbb R^{d-1}}I_1(x_0-x)c_d|x-z|^{-d}\,dv\notag\\
 &=\frac{c_d}{\pi b_0}\int_a^\infty
    \frac{h+b_0}{(|v_0|^2+(h+b_0)^2)^{d/2}}\,
                         \frac{dh}{\sqrt{h^2-a^2}}.\label{slic}
\end{align}
With $u=\sqrt{h^2-a^2}$ in \eqref{slic}, the numerator including the differential becomes
$(1+b_0/\sqrt{a^2+u^2})\,du$, and the denominator is at least
$(R_0^2+u^2)^{d/2}$. Now
$\int_0^\infty(1+u^2/R_0^2)^{-d/2}du\le CR_0/\sqrt d$.
For the term containing $(a^2+u^2)^{-1/2}$ we have two bounds. We can use $(a^2+u^2)^{-1/2}\le1/a$ throughout. Alternatively, after extracting $R_0^{-d}$, we can split the integral at $R_0$. The part below $R_0$ is at most $C(1+\log_+(R_0/a))$, and the tail is at most $\int_{R_0}^\infty(u/R_0)^{-d}\,du/u=1/d$. Recalling from \eqref{pgeom} that $b_0\asymp\ell p$ and $R_0\lesssim p$, we conclude that
\begin{equation}\label{psli}
 Q(t)\le C c_dR_0^{-d}
 \left[\frac1{\ell\sqrt d}+
       \min\left\{1+\log_+(Cp/a),\frac{Cp}{\sqrt d\,a}\right\}\right].
\end{equation}
Here $\log_+v=\max(0,\log v)$. The logarithmic bound makes the singularity at $a=|\tau-t|=0$ integrable in $t$.

\subsection{Heat annuli and completion of the proof}

We now integrate the slice estimate \eqref{psli} against $g$. This requires a comparison with the heat kernel on annuli of relative width $d^{-1/2}$. For $A_0>0$ and $d\ge1$,
\begin{equation}\label{pheat}
 H_{r_*^2/(2d)}(u_0)\asymp_{A_0} c_du_0^{-d}
 \qquad\text{if }|\log(u_0/r_*)|\le A_0/\sqrt d .
\end{equation}
Indeed, at $u_0=r_*$ the ratio of $c_du_0^{-d}$ to the heat kernel is $c_d(2\pi e/d)^{d/2}\asymp1$ by Stirling's formula. At $u_0=r_*e^v$ the ratio acquires the factor $\exp[d(e^{2v}-1-2v)/2]$, which lies between $1$ and a constant depending only on $A_0$ for $|v|\le A_0/\sqrt d$.

Consider $I_1*\nu_{y,p}(x_0)$ with $x_0\in F'\cap E_p(y)$. For each $z$ in the range $K_0\sqrt w<\ell\le5$, we bound the angular factor in \eqref{nudef} by $C\ell^2$, then enlarge the domain in $x$ to the slab $|t|\le U_0$. The coordinates of this slab and its width depend on $z$. By \eqref{slic} it suffices to prove
\begin{equation}\label{ppot}
 \int_{\{K_0\sqrt w<\ell\le5\}}g(z)\frac{\ell^2}{p}
       \int_{-U_0}^{U_0}Q(t)\,dt\,dz
       \le C w\big[1+\log(1/w)+L_d\big]L_d\Lambda.
\end{equation}
We first consider $r'\ge r/2$. Using $R_0\ge r'$ in \eqref{psli}, we claim that
\begin{equation}\label{pmix}
 \frac{\ell^2}{p}\int_{-U_0}^{U_0} Q(t)\,dt
 \le Cw c_d(r')^{-d}
 \begin{cases}
  1+\log(1/w),& |r'-r|\le A_1p/\sqrt d,\\
  d^{-1/2}+p/(\sqrt d\,|r'-r|),& |r'-r|>A_1p/\sqrt d,
 \end{cases}
\end{equation}
for a sufficiently large absolute constant $A_1$. After extracting $c_d(r')^{-d}$, the first term in \eqref{psli} contributes at most $Cw/\sqrt d$, since the interval has length $2U_0$. For the second term in the central range, the integral of $1+\log_+(Cp/|\tau-t|)$ is largest when the interval is centered at $\tau$. It is therefore bounded by $CU_0[1+\log_+(Cp/U_0)]$. Including the factor $\ell^2/p$, we obtain the bound $Cw\ell[1+\log_+(C\ell/w)]\le Cw[1+\log(1/w)]$.

Outside the central range, \eqref{pgeom} implies $|\tau|\ge c|r'-r|/\ell$. As $dw\le c$, we have $U_0\le |\tau|/2$ and hence $a\ge c|r'-r|/\ell$. The bound $Cp/(\sqrt d\,a)$ in \eqref{psli}, together with $\ell\le5$, now proves the second line of \eqref{pmix}.

To integrate \eqref{pmix} on the central interval, we use the single Gaussian $H_{r^2/(2d)}(r')$ from \eqref{pheat}. For the rest of the range $r/2\le r'\le7p$, split at $r$ and subdivide each side into annuli $jp/\sqrt d\le |r'-r|<(j+1)p/\sqrt d$, $1\le j\le C\sqrt d$. For each nonempty radial interval in $[r/2,7p]$, let $r_*$ be its midpoint. Then $r_*\asymp p$ and $|r'-r_*|\le p/\sqrt d$. By \eqref{pheat}, $c_d(r')^{-d}\le C H_{r_*^2/(2d)}(r')$ there. Integration against $g(z)$ and summation over these annuli yield the bound
\begin{equation*}
 Cw\left[1+\log(1/w)+\sum_{j\le C\sqrt d}(d^{-1/2}+j^{-1})\right] H_*g(x_0)
 \le Cw[1+\log(1/w)+L_d]H_*g(x_0).
\end{equation*}

It remains to consider $r'<r/2$. In this range $\ell\ge1/2$, $|\tau|\ge cp/\ell$, and $a\ge cp/\ell$. Moreover, $R_0\ge r-2U_0$ and $dU_0/p\le Cdw$, so $R_0^{-d}\le Cr^{-d}$. Integrating \eqref{psli} in $t$ with weight $\ell^2/p$ now yields at most $Cw c_d r^{-d}/\sqrt d$. By \eqref{pheat} at radius $r$ and the radial monotonicity of the Gaussian, $c_d r^{-d}\le C H_{r^2/(2d)}(r')$. The contribution of this range is thus at most $CwH_*g(x_0)$.

Since $H_*g(x_0)\le C L_d\Lambda$ on $F'$ by \eqref{bexc}, the two ranges together prove \eqref{ppot}. For the full density \eqref{nudef}, we add the bound $Cw\Lambda$ already obtained when $D\le K_0\sqrt w\,p$. We conclude that
\[
 I_1*\nu_{y,p}(x_0)\le A,\qquad
 A=Cw[1+\log(1/w)+L_d]L_d\Lambda
 \quad\text{for }\nu_{y,p}(x_0)\,dx_0\text{ almost every }x_0 .
\]
On $p\le s<2p$, we have $p c_ds^{-d}\le(d-1)I_1(x-y)$. A pair in this shell contributes to \eqref{near} only if $R<d^2P$, in which case $p\le s\le(1+d^2)P$. We may therefore integrate \eqref{ppair} over $F'$, sum the shell contributions, and apply \eqref{csurv} with $B=1+d^2$ to obtain
\begin{align*}
 \int_{F'} E_w^2
 &\le(d-1)\int_\Omega g(y)\sum_p (I_1*\nu_{y,p})(y)\,
       \pi_y\{z_*:P\ge p/(1+d^2)\}\,dy\\
 &\le Cd A L_d M .
\end{align*}
This proves \eqref{pthin} upon substitution of the value of $A$.

\begin{proof}[Completion of the proof of Theorem~\ref{main}]

Choose $w=c_0d^{-7/6}$ for a sufficiently small absolute constant $c_0>0$. The term involving $w$ in \eqref{bbest} is then $O(d^{-1/12}L_d^{1/2})$, and the square root of the coefficient of $\Lambda M$ in \eqref{pthin} is $O(d^{-1/12}L_d^{3/2})$. Both are bounded uniformly in $d$. It follows that
\begin{equation}\label{pfin}
 \|\mathbf1_{F'}\mathcal N(k_w*_u A)\|_2+\|\mathbf1_{F'}E_w\|_2
 \le C\sqrt{\Lambda M}.
\end{equation}
In the bound \eqref{recov} for $\mathcal B_dg$, we have $\|E_1\|_1\le CM$ by \eqref{pfar}, while \eqref{pfar} and \eqref{pfin} bound all other terms in $L^2(F')$ by $C\sqrt{\Lambda M}$. For the remaining part of $f=g+f\mathbf1_F$, we use \eqref{bLtwo} and $\|f\mathbf1_F\|_2^2\le\Lambda M$. Chebyshev's inequality and the estimate for $|(F')^c|$ in \eqref{bexc} now imply
\begin{equation*}
 |\{x:\mathcal B_df(x)>\lambda\}|
 \le \frac{CM}{\Lambda}+\frac{CM}{\lambda}
                  +\frac{C\Lambda M}{\lambda^2}.
\end{equation*}
Taking $\Lambda=\lambda$ proves \eqref{w11} for bounded nonnegative $f$ in sufficiently large dimension. For the remaining finitely many dimensions, apply the classical weak-type bound for each scalar maximal truncation (cf. e.g. \cite[Corollary~3.8.2]{Gra24}) and take the union of its exceptional sets at height $\lambda/\sqrt d$.

For arbitrary nonnegative $f\in L^1$, let $f_n=\min(f,n)$. At each $t>0$ the truncated kernel is bounded, so $R_{i,t}f_n\to R_{i,t}f$ uniformly in $x$. Fatou's lemma first yields \eqref{w11} for every finite family of positive rational radii. Exhausting these radii proves the full estimate. For complex $f$, apply it to the positive and negative parts of $\operatorname{Re}f$ and $\operatorname{Im}f$. Their $L^1$ norms sum to at most $\sqrt2\|f\|_1$, and by subadditivity at least one of the four maximal functions exceeds $\lambda/4$ wherever $\mathcal B_df>\lambda$.
\end{proof}

\begin{corollary}\label{mweak}
Let $\nu$ be a finite complex Borel measure on $\mathbb R^d$, and define its hard truncations by integration against $d\nu$. Then $\|\mathcal B_d\nu\|_{1,\infty}\le C|\nu|(\mathbb R^d)$ with the constant of Theorem~\ref{main}.
\end{corollary}

\begin{proof}
By Fubini's theorem, for almost every $x$ we have $|\nu|(\partial B(x,t))=0$ simultaneously for all rational $t>0$. At such an $x$, mollification of $\nu$ by a positive approximate identity recovers each rational truncation by dominated convergence, since the truncated kernel is bounded and continuous off the sphere. The mollified densities have $L^1$ norms at most $|\nu|(\mathbb R^d)$. As above, applying Theorem~\ref{main} on finite families of rational radii and then exhausting them proves the bound for their supremum. Finally, for every $x$ the integral over $|x-y|>t$ is right continuous in $t>0$ by dominated convergence. Thus the same bound holds for all radii, including those for which the sphere carries mass.
\end{proof}

\section{Sharp dependence on the Lebesgue exponent}\label{s8}

\subsection{Strong bounds for the maximal truncations}

One coordinate suffices for the lower bound throughout $1<p<\infty$. To obtain the upper bound with independent radii, we compare the hard truncations with heat and Poisson integrals for $p\ge3$ and then interpolate with the estimate \eqref{w11}.

\begin{proposition}\label{pgrow}
For all $d\ge1$ and $1<p<\infty$,
\begin{equation}\label{pnorm}
 \cot\!\left(\frac{\pi}{2p^*}\right)
 \le \|\mathcal A_d\|_{p\to p}\le\|\mathcal B_d\|_{p\to p}\le Cp^* ,
\end{equation}
where $C$ is absolute.
\end{proposition}

\begin{proof}
Letting the radius tend to zero in the first coordinate, we find that $|R_1f|\le\mathcal A_d f$ almost everywhere. For $d\ge2$, the scalar norm formula of Ba\~nuelos and Wang \cite{BW95} yields the lower bound in \eqref{pnorm}. For $d=1$, we use Pichorides' norm formula for the Hilbert transform given in \cite[p.~165, parts~(i) and~(iv), and Section~4]{Pic72}.

For the upper bound, suppose first that $p\ge3$ and $f$ is a real Schwartz function. We compare the hard cutoff in \eqref{rdef} with the heat tail $Q_{i,a}$ in \eqref{qheat} at $a=t^2/(4\alpha)$. With $v=\log(|z|/t)$ and $G_\alpha$ as in \eqref{gdist}, the error in the radial cutoff satisfies
\begin{equation*}
 \delta_d(v):=\bigl|\mathbf1_{\{v>0\}}-G_\alpha(\alpha e^{2v})\bigr|
 \le \exp[-\alpha(e^{2v}-1-2v)].
\end{equation*}
Indeed, $\int_0^\infty e^{\lambda u}\gamma_\alpha(u)\,du=(1-\lambda)^{-\alpha}$ for $\lambda<1$. Exponential Markov bounds with $\lambda=1-e^{-2v}$ apply to the upper tail when $v>0$ and the lower tail when $v<0$. In particular, $\delta_d(v)\le e^{-cd\min\{v^2,|v|\}}$.

In polar coordinates, apply Cauchy--Schwarz to the spherical integral. Since $\int\theta_i^2\,d\sigma_{d-1}(\theta)=1/d$, it follows that
\begin{equation*}
 |R_{i,t}f(x)-Q_{i,t^2/(4\alpha)}f(x)|
 \le C\int_{\mathbb R}\delta_d(v)\bigl(S_{te^v}|f|^2(x)\bigr)^{1/2}\,dv .
\end{equation*}
Let $\eta=d^{-1/2}$ and $I_k=[k\eta,(k+1)\eta)$. For $\tau_k=(te^{(k+1/2)\eta})^2/(2d)$, the radial heat estimate \eqref{pheat} implies that $|S^{d-1}|(te^v)^d H_{\tau_k}(te^v)\gtrsim\eta^{-1}$ on $I_k$. Integrating against this Gaussian, we find that
\begin{equation*}
 \int_{I_k}S_{te^v}|f|^2(x)\,dv\le C\eta H_*(|f|^2)(x).
\end{equation*}
For $|k|\ge2$ we have $\sup_{I_k}\delta_d\le C\exp[-c\min(k^2,\sqrt d\,|k|)]$. These bounds are summable in $k$ uniformly in $d$. Applying Cauchy--Schwarz on each $I_k$ and summing, we bound the error by $C d^{-1/2}(H_*(|f|^2)(x))^{1/2}$ in every coordinate, uniformly in $t$. The factor $d^{-1/2}$ cancels when we sum the squared coordinate maxima.

We next compare $Q_{i,a}$ with the conjugate Poisson integral. In the formula for $P_{\sqrt a}R_i$ in \eqref{cpker}, we subordinate the Poisson kernels by \eqref{subor} with $m=1$ and integrate in height. Thus
\begin{equation*}
 P_{\sqrt a}R_if=-\frac1{\sqrt\pi}\int_0^\infty
 e^{-a/(4u)}\partial_{x_i}H_u f\,\frac{du}{\sqrt u}.
\end{equation*}
By scaling, $\|\mathbf1_{\{u>a\}}-e^{-a/(4u)}\|_{L^2(du/u)}$ is independent of $a$ and finite. Cauchy--Schwarz in $du/u$ therefore bounds $|Q_{i,a}f-P_{\sqrt a}R_if|$ by $C(\int_0^\infty|\partial_{x_i}H_u f|^2\,du)^{1/2}$. Together with the radial cutoff estimate, this proves the pointwise bound
\begin{equation}\label{predu}
 \mathcal B_df\le \mathcal N(P_s\mathbf R f)+C G_Hf+
                        C\bigl(H_*(|f|^2)\bigr)^{1/2},
 \qquad G_Hf:=\left(\int_0^\infty|\nabla H_u f|^2\,du\right)^{1/2}.
\end{equation}
For $p\ge3$, the $L^p$ norm of the last term in \eqref{predu} is bounded by $C\|f\|_p$. Indeed, Jensen's inequality bounds $[H_*(|f|^2)]^{p/3}$ by $H_*(|f|^{2p/3})$, to which we apply the heat maximal theorem at the single exponent $3/2$ (see e.g. \cite[Chapter~III, Section~3]{SteinLP70}). For the square function and the Poisson maxima we use the estimates \eqref{kgoal} from Section~\ref{skern}. We conclude that $\|\mathcal B_df\|_p\le Cp\|f\|_p$ for real Schwartz functions, and subadditivity extends the estimate to complex Schwartz functions. Along a Schwartz approximation to $f$ in $L^p$, H\"older's inequality ensures convergence of the truncated integrals at each positive radius. Fatou's lemma, first for finite sets of radii and then along an exhaustion of the positive rationals, proves the estimate for every $f\in L^p$.

For $1<p\le3$ and $\lambda>0$, apply \eqref{w11} to $f\mathbf1_{\{|f|>\lambda\}}$ and the bound at $p=4$ to $f\mathbf1_{\{|f|\le\lambda\}}$. Integrating the resulting bound for $|\{\mathcal B_df>2\lambda\}|$ against $p\lambda^{p-1}\,d\lambda$, we obtain
\[
 \|\mathcal B_d f\|_p^p\le C2^pp\big[(p-1)^{-1}+(4-p)^{-1}\big]\|f\|_p^p .
\]
Taking the $p$th root proves the upper bound for $\mathcal B_d$ in \eqref{pnorm}.
\end{proof}

\subsection{Estimates for iterated heat and Poisson averages}\label{skern}

For the two extension terms in \eqref{predu}, we need different powers of $p$. We shall prove that, for every real Schwartz function $f$ and $p\ge3$,
\begin{equation}\label{kgoal}
 \|G_H f\|_p\le C\sqrt p\,\|f\|_p,\qquad
 \|\mathcal N(P_s\mathbf Rf)\|_p\le Cp\,\|f\|_p .
\end{equation}
Both estimates follow by iteration of kernel averages. Integrate the initial spatial variable against Lebesgue measure and each subsequent increment against its probability measure. We write $\langle\cdot\rangle_r$ for the average over the $(r+1)$st increment, conditional on all variables through step $r$.

We first estimate the accumulated squared increments. Suppose that $F_r$ depends only on the variables through step $r$ and that $\langle F_{r+1}\rangle_r=F_r$. Set
\begin{equation*}
 e_r:=\langle(F_{r+1}-F_r)^2\rangle_r,\qquad
 E_r:=\sum_{j<r}e_j .
\end{equation*}
For $q=p/2$, use $e_r=\langle F_{r+1}^2\rangle_r-F_r^2$ and the fact that $E_{r+1}$ depends only on the variables through step $r$. By Fubini's theorem and summation by parts,
\begin{align*}
 \int E_n^q
 &\le q\sum_{r<n}\int E_{r+1}^{q-1}e_r\\
 &=q\sum_{r<n}\int E_{r+1}^{q-1}(F_{r+1}^2-F_r^2)
 \le q\int E_n^{q-1}F_n^2.
\end{align*}
The last inequality follows from monotonicity of $E_{r+1}^{q-1}$ in $r$. H\"older's inequality now yields
\begin{equation}\label{fener}
 \|E_n^{1/2}\|_p\le \sqrt{p/2}\,\|F_n\|_p.
\end{equation}

To control coordinate maxima, suppose that each increment includes a sign $\varepsilon_{r+1}=\pm1$ with conditional probabilities $1/2$. For real vectors $a_r$ depending only on the variables through step $r$, write
\begin{equation*}
 N_j:=\sum_{r<j}\varepsilon_{r+1}a_r,\qquad
 V:=\left(\sum_i\max_{j\le n}|N_{j,i}|^2\right)^{1/2},
 \qquad S^2:=\sum_{r<n}|a_r|^2.
\end{equation*}
Then
\begin{equation}\label{smax}
 \|V\|_p\le C\sqrt p\,\|S\|_p.
\end{equation}
To see where the uniformity in the number of coordinates enters, put $m_{r,i}^{\pm}=\max_{j\le r}(\pm N_{j,i})$, with $N_0=0$, and $\Delta N_{r,i}=N_{r+1,i}-N_{r,i}$. A summation by parts at the times when $m_{r,i}^+$ increases yields
\begin{equation*}
 m_{n,i}^+N_{n,i}
 =\sum_{r<n}m_{r,i}^+\Delta N_{r,i}
  +\sum_{r<n}m_{r+1,i}^+(m_{r+1,i}^+-m_{r,i}^+)
 \ge\sum_{r<n}m_{r,i}^+\Delta N_{r,i}+\tfrac12(m_{n,i}^+)^2.
\end{equation*}
Apply the same estimate to $-N_{j,i}$ and use $2ab\le a^2/2+2b^2$. Expanding $N_{n,i}^2$ and summing over $i$, we obtain a bound in terms of a single signed sum,
\begin{equation*}
 V^2\le8S^2+\sum_{r<n}\varepsilon_{r+1}c_r,\qquad
 c_r:=\bigl(16N_r+4(m_r^--m_r^+)\bigr)\cdot a_r .
\end{equation*}
In particular, $(\sum_r c_r^2)^{1/2}\le20VS$.

Denote the partial sums of that scalar term by $Z_j=\sum_{r<j}\varepsilon_{r+1}c_r$, and set $Z^*=\max_{j\le n}|Z_j|$ and $D^2=\sum_{r<n}c_r^2$. After any initial segment, stop before adding a coefficient that would bring the sum of the subsequent squared coefficients above $u$, or when the absolute value of the subsequent signed sum first reaches $b$. Denote this stopped sum by $Z$ and the sum of its squared coefficients by $D_{\rm stop}^2$. Since $\cosh v\le e^{v^2/2}$, successive averaging shows that $\exp(\theta Z-\theta^2D_{\rm stop}^2/2)$ has conditional mean at most one. Taking $\theta=\pm b/u$, we bound the probability of reaching $b$ before exceeding $u$ by $2e^{-b^2/(2u)}$.

On $\{D\le\delta\lambda\}$, the first crossing of $\lambda$ overshoots by at most $\delta\lambda$. If $0<\delta\le1/2$, a subsequent crossing of $2\lambda$ therefore requires a further displacement of at least $\lambda/2$. Partition according to the first crossing and apply the preceding bound with $b=\lambda/2$ and $u=\delta^2\lambda^2$. We obtain
\[
 |\{Z^*>2\lambda\}|\le |\{D>\delta\lambda\}|
       +2e^{-1/(8\delta^2)}|\{Z^*>\lambda\}|.
\]
After integration against $q\lambda^{q-1}\,d\lambda$, the choice $\delta=[8(q+2)\log2]^{-1/2}$ allows us to absorb the last term. Hence $\|Z^*\|_q\le C\sqrt q\,\|D\|_q$ for $q\ge1$. Since $D\le20VS$, H\"older's inequality with $q=p/2$ implies
$\|V\|_p^2\le8\|S\|_p^2+C\sqrt p\,\|V\|_p\|S\|_p$. This proves \eqref{smax}. When $\|S\|_p<\infty$ the absorptions are justified by $V\le\sqrt n S$ and $Z^*\le\sqrt n D$.

We apply the energy estimate \eqref{fener} first to the heat square function. Given $T>0$, let $h=T/n$ and $t_r=T-rh$. Starting with Lebesgue measure in $x_0$, integrate each increment $x_{r+1}-x_r$ against $H_h$. For $F_r=H_{t_r}f(x_r)$, the semigroup property implies $\langle F_{r+1}\rangle_r=F_r$, and translation invariance yields $\|F_n\|_p=\|f\|_p$. At step $r$ the squared increment has conditional mean $e_r(x_r)$, where
\begin{equation*}
 e_r(x)=H_h(|H_{t_r-h}f|^2)(x)-|H_{t_r}f(x)|^2
       =2\int_0^h H_v(|\nabla H_{t_r-v} f|^2)(x)\,dv .
\end{equation*}
The integral identity follows by differentiating $H_v(|H_{t_r-v}f|^2)$. Now use $x_n$ in place of $x_0$ as the integration variable. By symmetry of the heat kernels, the average of $E_n$ with $x_n$ prescribed is $\sum_{r<n}H_{t_r}e_r(x_n)$. Jensen's inequality bounds its $L^{p/2}(\mathbb R^d)$ norm by $\|E_n\|_{p/2}$. Moreover,
\[
 |\nabla H_{2t_r}f|^2\le H_{t_r+v}(|\nabla H_{t_r-v}f|^2),\qquad 0\le v\le h .
\]
Thus $2h\sum_{r<n}|\nabla H_{2t_r}f|^2\le\sum_{r<n}H_{t_r}e_r$. As $n\to\infty$, the left side converges to $\int_0^{2T}|\nabla H_u f|^2\,du$. Letting $T\to\infty$ in the resulting $L^{p/2}$ estimate proves the heat bound in \eqref{kgoal}.

For the conjugate Poisson integrals, we use harmonic extension to the strip $-h<y<h$. Prescribe the same boundary function on both faces. At $y=0$ the extension is convolution with the kernel $\mathcal K_h$ whose Fourier multiplier is $\operatorname{sech}(h|\xi|)$. More generally, the multiplier at height $y$ is $\cosh(y|\xi|)/\cosh(h|\xi|)$. For smooth compactly supported boundary data, the extension tends to zero as $|x|\to\infty$, uniformly for $|y|\le h$ by Fourier inversion. The maximum principle therefore implies $\mathcal K_h\ge0$. At the origin, the multiplier and its first two derivatives determine the moments
\begin{equation*}
 \int \mathcal K_h=1,\qquad \int z_i \mathcal K_h(z)\,dz=0,\qquad
 \int z_i z_j \mathcal K_h(z)\,dz=h^2\mathbf1_{\{i=j\}}.
\end{equation*}
Since the multiplier is Schwartz, all these moments are absolutely convergent. We also have the identity
\begin{equation}\label{strip}
 \tfrac12 \mathcal K_h*(P_{(j+1)h}+P_{(j-1)h})=P_{jh},
\end{equation}
for integers $j\ge1$, with $P_0$ the identity.

Start at $(x_0,J)$ with Lebesgue measure in $x_0$ and an integer $J\ge1$. At each step average independently over a horizontal displacement $z$ of density $\mathcal K_h$ and a vertical sign $\varepsilon_{r+1}=\pm1$ with probabilities $1/2$. From $(x_r,j_r)$ with $j_r>0$ move to $(x_r+z,j_r+\varepsilon_{r+1})$. Once the path reaches level zero it stays there, with no further displacement. By \eqref{strip}, $U_j=P_{jh}f$ and $F_r=U_{j_r}(x_r)$ satisfy $\langle F_{r+1}\rangle_r=F_r$. For $j>0$ consider the spatial first moments
\begin{equation*}
 a_i(x,j):=\sum_{\varepsilon=\pm1}\frac12\int
     \frac{z_i}{h}\bigl(U_{j+\varepsilon}(x+z)-U_j(x)\bigr)\mathcal K_h(z)\,dz,
 \qquad a_i(x,0):=0 .
\end{equation*}
The moment identities make the functions $z_i/h$ orthonormal for the product of $\mathcal K_h(z)\,dz$ and the uniform sign measure. Bessel's inequality applied to $F_{r+1}-F_r$ yields $|a(x_r,j_r)|^2\le e_r$. With $a_r=a(x_r,j_r)$ in $N_j$, we therefore have $S^2\le E_n$. For each choice of increments, $x_n$ is a translate of $x_0$ and $P_{j_nh}$ is an $L^p$ contraction. Thus $\|F_n\|_p\le\|f\|_p$. Combining \eqref{fener} and \eqref{smax}, we conclude that $\|V\|_p\le Cp\|f\|_p$.

To recover the Poisson integrals at several heights, let $r_k$ be the first visit to level $k$, $1\le k<J$, or $n$ if the path has not visited $k$ by step $n$. We average $N_{r_k,i}$ against the terminal Poisson kernel. More precisely, for a Schwartz test function $g$ define $T^{(n)}_{i,k}f$ by
\begin{equation}\label{tstop}
 \int (T^{(n)}_{i,k}f)g=\int N_{r_k,i}\,(P_{j_nh}g)(x_n)\,d\mu.
\end{equation}
Here $d\mu$ is the product of Lebesgue measure in $x_0$ and the measures of the increments. Change variables from $x_0$ to $x_n$ and then to the integration variable of $g$ in $P_{j_nh}g$. For each value of this last variable the resulting average has mass one. Jensen's inequality for the norm $\ell_i^2(\ell_k^\infty)$ now bounds the $L^p$ norm of $(\sum_i\max_k|T^{(n)}_{i,k}f|^2)^{1/2}$ by $\|V\|_p$.

To identify $T^{(n)}_{i,k}$, expand $N_{r_k,i}$ in \eqref{tstop}. In the term from step $r$, iteration of \eqref{strip} replaces $P_{j_nh}g(x_n)$ by $P_{j_{r+1}h}g(x_{r+1})$ after averaging the subsequent steps. At level $j>k$ the signed average of this last function is $\mathcal K_h*(P_{(j+1)h}-P_{(j-1)h})g/2$, with multiplier $-e^{-jh\rho}\tanh(h\rho)$, where $\rho=|\xi|$. Differentiation of $\widehat{\mathcal K_h}$ yields $\widehat{a_i(\cdot,j)}=-e^{-jh\rho}\tanh(h\rho)\widehat{R_if}$. For each vertical path, translation removes the preceding horizontal increments from the integral in $x_0$. Thus the contribution from level $j$ at step $r$ has scalar factor $e^{-2jh\rho}\tanh^2(h\rho)$ multiplying the Riesz multiplier $-i\xi_i/\rho$, with weight equal to the probability of visiting $j$ at step $r$ before the first visit to $k$.

For $u=h\rho>0$, denote the sum of these scalar factors, with their weights, by $b_n(j)$ when the path starts at height $j\ge k$. Conditioning on the first vertical step yields $b_0=0$, $b_n(k)=0$, and
\begin{equation}\label{disc}
 b_{n+1}(j)=\tfrac12[b_n(j-1)+b_n(j+1)]+e^{-2ju}\tanh^2 u,\qquad j>k.
\end{equation}
Using $\cosh(2u)-1=2\sinh^2u$, one checks by substitution that
\[
 b(j)=\tfrac12\operatorname{sech}^2u\,(e^{-2ku}-e^{-2ju})
\]
is a bounded stationary solution of \eqref{disc}. Positivity of the averaging operator implies $0\le b_n\uparrow b_\infty\le b\le1/2$. The bounded sequence $b-b_\infty$ has zero second difference, so it is affine and hence constant on $\{k,k+1,\ldots\}$. It vanishes at $k$, which proves $b_\infty=b$. Dominated convergence in the Fourier integral for $f$ now yields the pointwise limit
\[
 T^{(n)}_{i,k}f\longrightarrow\tfrac12\operatorname{sech}^2(h|D|)(P_{2kh}-P_{2Jh})R_i f .
\]
Let $J\to\infty$. For any finite set of positive heights $s$, choose integers $k=k(h,s)\ge1$ such that $2kh\to s$ as $h\downarrow0$. Dominated convergence in the Fourier integrals and Fatou's lemma preserve the bound for the coordinatewise maxima over these heights. Exhausting the positive rationals and using continuity in $s$ proves the Poisson estimate in \eqref{kgoal}.

\section*{Acknowledgements}
\addcontentsline{toc}{section}{Acknowledgements}
S.-C. Mao is supported by the China Postdoctoral Science Foundation under Grant No.~2026M793367.

\bigskip
\noindent\textbf{Competing interests.}
The author declares no competing interests.

\bigskip
\noindent\textbf{Data availability.}
No datasets were generated or analyzed in this study.

\bigskip
\noindent\textbf{AI declaration.}
The author used AI tools to assist with the preliminary investigation and manuscript revision.

\par\bigskip
\pdfbookmark[1]{References}{refs}
\begingroup\linespread{1.00}\selectfont

\endgroup

\bigskip\bigskip
\noindent\begin{minipage}{\textwidth}

\bigskip
\noindent{Sheng-Chen Mao}\\
School of Mathematics and Statistics, Lanzhou University\\
No.~222 Tianshui South Road, Lanzhou 730000, P.~R. China\\
\textit{Emails:} \href{mailto:maoshengchen@lzu.edu.cn}{maoshengchen@lzu.edu.cn};
\href{mailto:maosci@163.com}{maosci@163.com}

\end{minipage}

\end{document}